\documentclass[a4paper,12pt]{elsarticle}
\journal{ArXiv}
\usepackage{amsmath}
\usepackage[english]{babel}
\usepackage{graphicx}
\usepackage[usenames,dvipsnames]{color} 
\usepackage[numbers]{natbib}
\usepackage[footnotesize,bf,sf,center]{caption}
\usepackage{float}
\usepackage[dvips]{epsfig}
\usepackage[bottom=2.5cm,top=2.5cm,left=3cm,right=2cm]{geometry}
\usepackage{calc}
\usepackage{soul}
\usepackage{amsmath,amsfonts,amssymb}
\usepackage{subfigure}
\usepackage{setspace}
\usepackage{amssymb}
\usepackage{mathrsfs}
\usepackage{algorithm}
\usepackage{algpseudocode}
\usepackage{dsfont}
\usepackage{epstopdf}
\usepackage{subfig}%
\usepackage[utf8]{inputenc}
\usepackage{nicefrac,xfrac}
\newcommand{\ve}[1]{\mbox{\boldmath $#1$}}
\newtheorem{theorem}{Theorem}[section]

\newtheorem{lemma}[theorem]{Lemma}

\newtheorem{remark}{Remark}[section]

\newdefinition{rmk}{Remark}
\newcommand{\proof} [1]{ \noindent {\bf Proof.} #1 \hfill\rule{0.5em}{1.2ex} \par\medskip}
\usepackage{graphicx}
\begin{document}
% -------------------------------------------------------------------------------------------
\begin{frontmatter}

  \renewcommand\arraystretch{1.0}
      %\title{\textbf{Implicit-explicit discretisations of the viscous term in generalised  incompressible flows with variable viscosity}}

    \title{\textbf{Unconditionally stable, linearised IMEX schemes for incompressible flows with variable density}}
      
    \author{
   {\bf Nicolás Espinoza-Contreras} $^{1,2,3}$, \ 
   {\bf Gabriel R.~Barrenechea} $^4$,\
   {\bf Ernesto Castillo} $^5$, 
        and \  
   {\bf Douglas R.~Q.~Pacheco} $^{1,2,3,}$
   \footnote[0]{{\sf Email address:} {\tt pacheco@ssd.rwth-aachen.de}, 
     {\sf corresponding author}}\\
   {\small ${}^{1}$ Chair for Computational Analysis of Technical Systems, RWTH Aachen University, Germany}\\
   {\small ${}^{2}$ Chair of Methods for Model-based Development in Computational Engineering, RWTH}\\
   {\small ${}^{3}$ Center for Simulation and Data Science (JARA-CSD), RWTH Aachen University, Germany}\\
    {\small ${}^{4}$ Department of Mathematics and Statistics, University of Strathclyde, Glasgow, Scotland} \\
    {\small ${}^{5}$ Department of Mechanical Engineering, University of Santiago de Chile, Santiago, Chile}}
\begin{abstract}
%In inhomogeneous fluid flows, even if incompressible, viscosity and density can vary in space and time. \textcolor{blue}{As a consequence, the diffusive term can no longer be written as a simple vector Laplacian, but the full symmetric gradient  form needs to be kept. This fact leads to a coupling of all the components of the velocity, which is a computational drawback.} While a fully explicit temporal discretisation of the viscous term is usually prohibitive, \textcolor{blue}{for a problem with variable viscosity (but constant density)} we have recently shown that implicit-explicit (IMEX) treatments can lead to unconditionally stable methods. So,

For the incompressible Navier--Stokes system with variable density and viscosity, we propose and analyse an IMEX framework treating the convective and diffusive terms semi-implicitly. This extends to variable density and second order in time some methods previously analysed for variable viscosity and constant density. We present three new schemes, both monolithic and fractional-step. All of them share the methodological novelty that the viscous term is treated in an implicit-explicit (IMEX) fashion, which allows decoupling the velocity components. Unconditional temporal stability is proved for all three variants. Furthermore, the system to solve at each time step is linear, thus avoiding the costly solution of nonlinear problems even if the viscosity follows a non-Newtonian rheological law. Our presentation is restricted to the semi-discrete case, only considering the time discretisation. In this way, the results herein can be applied to any spatial discretisation. We validate our theory through numerical experiments considering finite element methods in space. The tests range from simple manufactured solutions to complex two-phase viscoplastic flows.

%In inhomogeneous fluid flows, even if incompressible, viscosity and density can vary in space and time. Therefore, the diffusive term cannot be written as a simple vector Laplacian, which brings computational disadvantages. While a fully explicit temporal discretisation of the viscous term is usually prohibitive, we have recently shown that implicit-explicit (IMEX)   treatments can lead to unconditionally stable methods. In this article, we propose and analyse an IMEX framework that can decouple the computation of density, velocity and pressure, as well as of the individual velocity components. Here, recent results obtained for first-order stepping and constant density are now extended to the case of variable density and temporal order up to two. We present three new schemes for the incompressible Navier--Stokes system with variable density, including an IMEX fractional-step method. Unconditional temporal stability is proved for all three variants. Furthermore, the system to solve at each time step is fully linearised; that is, there is no need to use Newton or Picard iterations, even if the viscosity follows complex rheological laws. With respect to state-of-the-art numerical methods for variable-density incompressible flows, the key methodological novelty here is an IMEX treatment of the (variable-viscosity) diffusive term. Our numerical tests use finite elements in space and range from simple manufactured solutions to complex two-phase viscoplastic flows.
\end{abstract}
\end{frontmatter}

\section{Introduction}
In incompressible flows, density variability occurs mainly in the presence of different fluid phases, be they immiscible or not. Since the seminal work by \citet{Guermond2000}, many numerical methods have been developed for the incompressible variable-density Navier--Stokes equations. For the temporal discretisation, the predominant approach is to use extrapolation to decouple the density and momentum equations. Most articles on numerical analysis for problems with variable density \cite{Guermond2000,Guermond2011,An2020,Chen2020,Cai2021} assume constant dynamic viscosity. This simplifying assumption facilitates the analysis since it makes it simpler to discretise the viscous term in a fully implicit manner. However, since the dynamic viscosity is proportional to the density, and, of course, since different fluids usually have different kinematic viscosities, variable density implies variable viscosity. So, in this work we shall also address variable viscosity and how to handle it numerically to preserve the good algorithmic properties observed for homogeneous flows.

%Interestingly, what most articles focusing on the corresponding numerical analysis neglect \cite{Guermond2000,Guermond2011,An2020,Chen2020,Cai2021} is that variable density implies variable viscosity---not only because different fluids usually have different kinematic viscosities, but also because the dynamic viscosity is proportional to the density. So, by assuming constant dynamic viscosity ``for simplicity'' of analysis, most theoretical works treat the viscous term implicitly. In the present work, on the other hand, we shall address variable viscosity and how to handle it numerically to preserve good algorithmic properties observed for homogeneous flows.  

For problems with variable density, decoupling strategies are usually preferred over monolithic discretisations.  Density-momentum decoupling is the standard in variable-density solvers, with pressure segregation also being a very common approach. As a matter of fact, there is a vast body of literature on corresponding fractional-step schemes, and we can sort most of these methods into three families. The first are incremental pressure-correction schemes \cite{Guermond2000, Deteix2022}, which are variable-density extensions of classical projection methods. The second family are Gauge--Uzawa methods \cite{Chen2020,Pyo2007,Liu2015}, which also belong to the class of projection methods. The third popular approach are penalty methods, which perturb the incompressibility constraint using a pressure increment Laplacian \cite{Guermond2011,Guermond2009,Salgado2013,Axelsson2015}. An advantage of the latter is that there is no need to solve a variable-coefficient pressure equation. Furthermore, no second-order extensions have been proven stable for those two families of projection methods, although numerical evidence often indicates good stability \cite{Guermond2000,Pyo2007,Wu2017}. A less common, yet very accurate alternative, are consistent splitting schemes, which are based on a fully consistent pressure equation and therefore eliminate splitting errors and numerical boundary layers \cite{Li2021,Pacheco2022CompMech}. In the variable-density context, pressure-correction methods usually require solving two pressure Poisson problems per time step \cite{Guermond2000}. An interesting variant recently proposed by \citet{Deteix2022}, from which we take inspiration herein, requires solving only one pressure Poisson problem per time step.

In this work, our main interest is to propose and analyse IMEX schemes. In CFD, the term IMEX usually refers to temporal discretisations that make convection (semi-)explicit while keeping viscous terms implicit, see e.g.~the analyses and methods presented very recently by Burman and co-workes \cite{Burman2022,Burman2023, Burman2024}. For variable density, too, most of the literature treats the convective term semi-implicitly to linearise the momentum equation and decouple the density update from the velocity-pressure step \cite{Guermond2000}. Now, when density is not constant, neither is viscosity, and in this case the diffusive term cannot be written as a simple vector Laplacian, which means that an implicit treatment will end up coupling all the velocity components, even when using a pressure segregation method. These considerations motivate the search for schemes that decouple the velocity components and linearise the systems to be solved at each time step. Thus, IMEX discretisations that make the transpose velocity gradient or other coupling terms explicit have both theoretical and practical relevance. This was the main motivation for our recent work \cite{Barrenechea2024}, where first-order IMEX schemes built under this premise were analysed for problems with variable viscosity (but constant density). Interestingly, not many papers to date are devoted to IMEX schemes for variable-viscosity problems \cite{Stiller2020,Guesmi2023,Yakoubi2024}. The present work generalises our recent approach to the variable-density case and also presents a second-order extension. We propose three different IMEX methods: a first- and a second-order method decoupling only density transport and momentum, and a first-order projection method fully decoupling density, pressure and velocity (components). We shall prove that all of our schemes, which are based on a consistent reformulation of the viscous term, are unconditionally stable in time.

The remaining content in this article is presented as follows. Section \ref{sec_pre} introduces the model problem and underlying assumptions, as well as useful notation, identities and inequalities. Section \ref{sec_schemes} presents our new IMEX schemes, whose temporal stability is then analysed in Section \ref{sec_analysis}. Section \ref{sec_impl} briefly addresses implementation matters, followed by numerical examples in Section \ref{sec_examples}. We finally draw concluding remarks in Section \ref{sec_Conclusion}.

\section{Preliminaries}\label{sec_pre}
\subsection{Model problem}
For a finite time interval $(0,T]$ and a domain $\Omega\subset\mathbb{R}^{d}$, $d=2$ or $3$, the variable-density incompressible Navier--Stokes equations can be written as
\begin{align}
\partial_t\rho + \ve{u}\cdot\nabla\rho  &= 0 && \text{in} \ \ \Omega\times(0,T]=:Q\, ,\label{mass}\\
\rho\partial_t\ve{u} + (\rho\nabla\ve{u})\ve{u} - \nabla\cdot(2\rho\nu\nabla^{\mathrm{s}}\ve{u}) + \nabla p &= \ve{f} && \text{in} \ \ \Omega\times(0,T]\, ,\label{momentum}\\
\nabla\cdot\ve{u} &= 0  && \text{in} \ \ \Omega\times(0,T]\, , \label{incompressibility}
\end{align}
%\ve{u} &= \ve{0} && \text{on} \ \ \Gamma\times(0,T]\, ,\label{DirichletBC}\\
%\ve{u} &= \ve{u}_0 && \text{at} \ \ t=0\, ,\\
%\rho &= \rho_0 && \text{at} \ \ t=0\, ,
where $\ve{u}$ is the velocity, $p$ is the pressure, $\rho>0$ is the density, $\nu>0$ is the kinematic viscosity, and $\ve{f}$ is a body force. Throughout this work, bold characters will refer to vectors. The dynamic viscosity $\mu:=\rho\nu$ is assumed to satisfy
\begin{equation}
    \mu \geq  \mu_{\mathrm{min}} > 0 \ \ \text{in} \ \bar{Q}\, ,
\end{equation}
where $\mu_{\mathrm{min}}$ is a constant usually determined by fluid properties (density and rheology). Initial and boundary conditions are also needed, which depends on the application, see e.g.~Ref.~\cite{Pacheco2022CompMech}. For the theoretical analysis, we will consider $\ve{u} = \ve{0}$ on $\partial\Omega$, so no boundary conditions are needed for $\rho$. In that case, if the initial density field $\rho_0(\ve{x})$ is such that $0 < \rho_{\text{min}} \leq \rho_0(\ve{x}) \leq \rho_{\text{max}}$ for all $\ve{x}\in\Omega$, then  
\begin{equation}
    \rho_{\text{min}} \leq \rho(\ve{x},t) \leq \rho_{\text{max}} \ \ \text{in} \ \bar{Q}\, ,
\end{equation}
as shown by \citet{Guermond2000}.

\subsection{Consistent reformulation}
When $\nu$ and $\rho$ can vary in space, the viscous term cannot be simplified to the Laplacian form $\mu\Delta\ve{u}$. On the other hand, the natural alternative, namely the stress-divergence form
\begin{equation}
\nabla\cdot\left(2\mu\nabla^{\mathrm{s}}\ve{u}\right) = \nabla\cdot(\mu\nabla\ve{u}) + \nabla\cdot(\mu\nabla^{\top}\ve{u}) \, ,
\label{SD}
\end{equation}
has computational disadvantages. In particular, the transpose gradient $\nabla^{\top}$ couples the velocity components $u^j$ because the $i$th entry in the vector $\nabla\cdot(\mu\nabla^{\top}\ve{u})$ involves all $u^j$ components, $j=1,\dots,d$. This does not happen for the other part of the viscous term, since the $i$th entry in $\nabla\cdot(\mu\nabla\ve{u})$ is simply
$\nabla\cdot(\mu\nabla u^i)$ -- which is why the Laplacian form is often preferred for dicretisation (when $\mu$ is constant) \cite{John2016}. In principle, the coupling could be avoided by treating the transpose gradient explicitly, but we have recently showed that this is generally unstable \cite{Barrenechea2024}. To overcome that, we use a consistent reformulation:
\begin{equation}
\begin{split}
\nabla\cdot\left(2\mu\nabla^{\mathrm{s}}\ve{u}\right) &\equiv \nabla\cdot\left(\mu\nabla\ve{u}\right) + \nabla^{\top}\ve{u}\nabla\mu+ \mu\nabla(\nabla\cdot\ve{u})  \\
    &= \nabla\cdot\left(\mu\nabla\ve{u}\right) + \nabla^{\top}\ve{u}\nabla\mu\, ,
\label{genLapStrong}    
\end{split}
\end{equation}
which we call generalised Laplacian form. The term $\nabla\cdot\left(\mu\nabla\ve{u}\right)$ is a variable-coefficient Laplacian, which at the fully discrete level leads to block-diagonal velocity matrices \cite{John2016,Barrenechea2024}. In this context, our idea is to treat explicitly $\nabla^{\top}\ve{u}$ and possibly $\mu$ (if the viscosity follows some non-Newtonian law), analysing also the implications of such an IMEX approach. To use this reformulation, we assume that  $\nabla\mu\in [L^\infty(\bar{Q})]^d$, which is essentially a Lipschitz condition on $\mu$.

For improved numerical stability, we will also re-write the convective terms consistently, as proposed by \citet{Guermond2000}. The reformulated equations are
\begin{align}
\partial_t\rho + \ve{u}\cdot\nabla\rho + \frac{\nabla\cdot\ve{u}}{2}\rho &= 0 \, ,\label{mass2}\\
\sqrt{\rho}\,\partial_t(\sqrt{\rho}\ve{u}) + (\rho\nabla\ve{u})\ve{u} + \frac{1}{2}[\nabla\cdot(\rho\ve{u})]\ve{u} - \nabla\cdot(\rho\nu\nabla\ve{u}) - \nabla^{\top}\ve{u}\nabla(\rho\nu) + \nabla p &= \ve{f} \, ,\label{momentum2}\\
\nabla\cdot\ve{u} &= 0 \, , \label{incompressibility2}
\end{align}
which are equivalent to \eqref{mass}--\eqref{incompressibility}.

\subsection{Useful notation, identities and inequalities}
We consider usual notation for Hilbert and Lebesgue spaces \cite{EG21-I}. The inner product and duality pairing in $L^2(\Omega)$ are denoted by $( \cdot ,\cdot )$ and $\langle \cdot ,\cdot \rangle$, respectively; we make no distinction between the product of scalar, vector- or tensor-valued functions. Moreover, we denote by $\| \cdot \|$, $\| \cdot \|_{-1}$ and $\| \cdot \|_{\infty}$ the $L^2(\Omega)$, $H^{-1}(\Omega)$ and $L^{\infty}(\bar{Q})$ norms, respectively. We assume that $\ve{f}\in [H^{-1}(\Omega)]^d$ for all $t$, and that the initial velocity condition satisfies $\ve{u}_0\in [H^1_0(\Omega)]^d$. To shorten notation, we define $\ve{U}:=\sqrt{\rho}\,\ve{u}$.

In our analysis, approximate or discrete-in-time values of the quantities at the different time steps will be denoted with sub-indices: $\ve{u}_{n}$, for instance, denotes the velocity approximation at the $n$th time step. The following identities will be useful:
\begin{align}
    2(\ve{v}_{n+1}-\ve{v}_{n},\ve{v}_{n+1}) &= \|\ve{v}_{n+1}\|^2 - \|\ve{v}_{n}\|^2 + \|\delta\ve{v}_{n+1}\|^2 , \label{identityBDF1}\\
     2(3\ve{v}_{n+1}-4\ve{v}_{n}+\ve{v}_{n-1},\ve{v}_{n+1}) &= \|\ve{v}_{n+1}\|^2 - \|\ve{v}_{n}\|^2 + \|\ve{v}_{n+1}^{\star}\|^2 - \|\ve{v}_{n}^{\star}\|^2 +\|\delta^2\ve{v}_{n+1}\|^2  ,
    \label{identityBDF2}
\end{align}
which holds for any scalar or vector-valued $\ve{v}_{n+1},\ve{v}_{n},\ve{v}_{n-1}$, where
\begin{align*}
\ve{v}_{n+1}^{\star} &:= 2\ve{v}_{n+1}-\ve{v}_{n}\, ,\\
\delta\ve{v}_{n+1}&:=\ve{v}_{n+1}-\ve{v}_{n}\, ,\\
\delta^2\ve{v}_{n+1}&:=\ve{v}_{n+1}-2\ve{v}_{n}+\ve{v}_{n-1} = \ve{v}_{n+1} - \ve{v}_{n}^{\star} \, .
\end{align*}
Another important fact is the skew-symmetry of the convective terms. For any $(r,\ve{v},\ve{w})$ sufficiently smooth, there holds
\begin{align}
\int_{\Omega} \Big[\ve{v}\cdot\nabla r + \frac{\nabla\cdot\ve{v}}{2}r\Big]r\,\mathrm{d}\Omega &= 0\, ,
\label{skewSymRho} \\
\int_{\Omega} \Big[(r\nabla\ve{v})\ve{w} + \frac{\nabla\cdot(r\ve{w})}{2}\ve{v}\Big]\cdot\ve{v}\,\mathrm{d}\Omega &= 0\, ,
\label{skewSymU}
\end{align}
provided that $\ve{v}\cdot\ve{n}=0$ on $\partial\Omega$ \cite{Guermond2000}.

We will analyse both coupled and fractional-step schemes. The stability analysis will rely on the following Gronwall inequality, proved by \citet{Heywood1990}.
\begin{lemma}[Discrete Gronwall inequality]\label{Lem:Gronwall-unconditional}
Let $N\in\mathbb{N}$, and $\alpha,B,a_{n},b_{n},c_{n}$ be non-negative numbers for $n=1,\ldots,N$. Let us suppose that these numbers satisfy
\begin{align}\label{N-1}
    a_{N} + \sum_{n=1}^{N}b_n \leq B + \sum_{n=1}^{N}c_n + \alpha\sum_{n=1}^{N-1}a_n  \, .
\end{align}
Then, the following inequality holds:
\begin{align}
    a_{N} + \sum_{n=1}^{N}b_n \leq \mathrm{e}^{\alpha N}\left( B + \sum_{n=1}^{N}c_n\right) \ \ \text{for} \ \ N\geq 1 \, .
\end{align}
\end{lemma}

IMEX stepping schemes can be constructed by using backward differentiation formulas (BDFs) and extrapolation rules of matching order. We will consider first- and second-order schemes, for which we have
\begin{align*}
    \text{$1^{\text{st}}$ order:}\ \ \partial_t(\sqrt{\rho}\, \ve{u})|_{t=t_{n+1}} &\approx \frac{1}{\tau}(\sqrt{\rho_{n+1}}\,\ve{u}_{n+1}-\sqrt{\rho_{n}}\,\ve{u}_{n}) \\
    \ve{u}_{n+1} &\approx \ve{u}_{n} \, , \\
    \text{$2^{\text{nd}}$ order:}\ \ \partial_t(\sqrt{\rho}\,\ve{u})|_{t=t_{n+1}} &\approx \frac{1}{2\tau}(3\sqrt{\rho_{n+1}}\,\ve{u}_{n+1}-4\sqrt{\rho_{n}}\,\ve{u}_{n} + \sqrt{\rho_{n-1}}\,\ve{u}_{n-1})\\
    \ve{u}_{n+1} &\approx 2\ve{u}_{n} - \ve{u}_{n-1} = \ve{u}_n^{\star} \, , 
\end{align*}
where $\tau = T/N >0$ denotes a constant time-step size, for simplicity. Similarly, the density equation \eqref{mass2} will be solved for $\rho_{n+1}$ either via
\begin{equation}
\frac{\rho_{n+1}-\rho_n}{\tau} + \ve{u}_n\cdot\nabla\rho_{n+1} + \frac{\nabla\cdot\ve{u}_{n}}{2}\rho_{n+1} = 0\, ,\label{densityBDF1}
\end{equation}
or through the second-order scheme
\begin{equation}
\frac{3\rho_{n+1}-4\rho_n+\rho_{n-1}}{2\tau} + \ve{u}^{\star}_n\cdot\nabla\rho_{n+1} + \frac{\nabla\cdot\ve{u}^{\star}_{n}}{2}\rho_{n+1} = 0\, .\label{densityBDF2}
\end{equation}

We will assume that the numerical density $\rho_n$ satisfies the lower bound
\begin{align}\label{Assumption-positivity}
    \rho_n \geq \varrho_{\mathrm{min}}  \ \text{a.e.~in}\ \Omega\, , \forall n={0,...,N},
\end{align}
where $\varrho_{\mathrm{min}}$ is a constant in $(0,\rho_{\mathrm{min}}]$.

\begin{remark} 
As stated by Guermond and Salgado \cite{Guermond2011,Guermond2009}, Assumption \eqref{Assumption-positivity} can be proven in the context of the semi-discrete analysis we carry out in this work. We start by defining the positive and negative parts of a function $v\in H^1(\Omega)$, a.e.~in $\Omega$, by
\begin{equation}
    v^+(\boldsymbol{x})=\max\{ v(\boldsymbol{x}),0\}\quad\textrm{and}\quad v^-(\boldsymbol{x})=v(\boldsymbol{x})- v^+(\boldsymbol{x})\, .
\end{equation}
Then, if we assume that $\rho_n^{}\ge\rho_{\mathrm{min}}$, 
multiplying \eqref{densityBDF1} by $(\rho_{n+1}-\rho_{\mathrm{min}})^-$, integrating over $\Omega$, using the fact that $\ve{u}_n$ is solenoidal and vanishes on the boundary (thus the convective term is skew-symmetric) and that $v^+$ and $v^-$ have complementary supports in $\Omega$, we arrive at the following:
\begin{equation}
    \int_\Omega \frac{[(\rho_{n+1}-\rho_{\mathrm{min}})^-]^2}{\tau} = \int_\Omega \frac{(\rho_{n}-\rho_{\mathrm{min}})(\rho_{n+1}-\rho_{\mathrm{min}})^-}{\tau}\le 0\,,
\end{equation}
and then \eqref{Assumption-positivity} follows with $\varrho_{\mathrm{min}}=\rho_{\mathrm{min}}$. Up to our best knowledge, the same proof cannot be repeated for the second-order scheme based on \eqref{densityBDF2}. In addition, in the fully discrete case, when the discrete velocity is not guaranteed to be divergence-free, the proof sketched above does not hold any longer. In both these situations, possible ways of deriving a discretisation that guarantees \eqref{Assumption-positivity} include methods either related to the FCT scheme, such as \cite{Bur15}, or a penalisation such as the one proposed in \cite{BurErn17}, both for the transport equation. Alternatively, the recent work \cite{ABP24} can also be used in this situation, where the restriction \eqref{Assumption-positivity} would be hardwired into the method regardless of the space or time discretisation.
\end{remark}

\section{Iteration-free IMEX methods}\label{sec_schemes}
In this article, we propose three different IMEX schemes. For the first two methods, velocity and pressure are computed simultaneously at each time step, after the density is updated through either \eqref{densityBDF1} or \eqref{densityBDF2}. Whenever the term ``coupled scheme`` is used herein, we mean the coupling between pressure and velocity; the density, on the other hand, will always be computed separately (beforehand) at each time step. Another common feature of all the methods presented herein is being fully linearised, that is, each time step requires the solution of linear subproblems only. 

\subsection{A first-order coupled scheme}
To derive a first-order linearised scheme, we use BDF1 for the time derivatives and extrapolate selected terms. The density is first updated via \eqref{densityBDF1}, followed by the velocity-pressure step:
\begin{multline}
\frac{\sqrt{\rho_{n+1}}}{\tau}(\sqrt{\rho_{n+1}}\,\ve{u}_{n+1}-\sqrt{\rho_n}\,\ve{u}_{n}) + (\rho_{n+1}\nabla\ve{u}_{n+1})\ve{u}_n + \frac{1}{2}[\nabla\cdot(\rho_{n+1}\ve{u}_n)]\ve{u}_{n+1}  \\
 -\nabla\cdot\left(\mu_{n+1}\nabla\ve{u}_{n+1}\right) 
 + \nabla p_{n+1} = \nabla^{\top}\ve{u}_{n}\nabla\mu_{n+1} +  \ve{f}_{n+1}\, ,\label{coupledBDF1}
 \end{multline}
\begin{equation}
 \nabla\cdot\ve{u}_{n+1} = 0\, .\label{divuBDF1}
\end{equation}
\begin{remark}
This scheme assumes that $\mu_{n+1}$ is already known at the point of computing $\ve{u}_{n+1}$ and $p_{n+1}$. This is often the case since the viscosity field is normally propagated with $\rho$, that is, $\mu_{n+1} = f(\rho_{n+1})$. In non-Newtonian applications, however, the viscosity may further depend locally on $\nabla^{\mathrm{s}}\ve{u}$ or $p$. In that case, we can avoid the nonlinearity by simply replacing $\mu_{n+1}$ with $\mu_n$ in Eq.~\eqref{coupledBDF1}, which does not affect the unconditional stability of the scheme. A more detailed discussion of the non-Newtonian case is given in Section \ref{sec_impl}.
\end{remark}
%\begin{remark}
%Using $\mu_n = \rho_n\nu_n$ instead of $\mu_{n+1} = \rho_{n+1}\nu_{n+1}$ may seem unnecessary, since $\rho_{n+1}$ is already known at the point of updating pressure and velocity. Although it is indeed possible, in many cases, to use $\mu_{n+1}$, the scheme presented above allows for generalised Newtonian rheologies, where the viscosity $\nu$ can be a function of $\nabla^{\mathrm{s}}\ve{u}$ and/or $p$. In that case, using $\nu_n$ circumvents that nonlinearity. Either way, the extrapolation $\nu_{n+1}\approx \nu_n $ is consistent with the order $\mathcal{O}(\tau)$ of the scheme. 
%end{remark}

\subsection{A second-order coupled scheme}
To extend the method to order two in time, we switch to BDF2 and second-order extrapolations. After computing $\rho_{n+1}$ through \eqref{densityBDF2}, we update velocity and pressure via
\begin{multline}
\frac{\sqrt{\rho_{n+1}}}{2\tau}(3\sqrt{\rho_{n+1}}\,\ve{u}_{n+1}-4\sqrt{\rho_{n}}\,\ve{u}_{n} + \sqrt{\rho_{n-1}}\,\ve{u}_{n-1}) + (\rho_{n+1}\nabla\ve{u}_{n+1})\ve{u}_n^{\star} + \frac{1}{2}[\nabla\cdot(\rho_{n+1}\ve{u}_n^{\star})]\ve{u}_{n+1}  \\
 -\nabla\cdot\left(\mu_{n+1}\nabla\ve{u}_{n+1}\right) 
 + \nabla p_{n+1} = \nabla^{\top}\ve{u}_{n}^{\star}\nabla\mu_{n+1} +  \ve{f}_{n+1}\, , \label{coupledBDF2}
 \end{multline}
\begin{equation}
 \nabla\cdot\ve{u}_{n+1} = 0\, . \label{divuBDF2}
\end{equation}
%where $(\nabla^{\top}\ve{u}_{n}\nabla\;.mu_{n})^{\star} := 2\nabla^{\top}\ve{u}_{n}\nabla\mu_{n} - \nabla^{\top}\ve{u}_{n-1}\nabla\mu_{n-1}$. 
As in the first-order case, $\mu_{n+1}$ can be replaced by a suitable extrapolation -- now second-order -- if the viscosity follows a nonlinear rheology (see Section \ref{sec_impl}).
\begin{remark}
    Since this is a multi-step method, we use a single first-order step as initialisation to compute $\rho_1,\ve{u}_1,p_1$. This is a standard approach and does not affect the (global-in-time) second order of the scheme. 
\end{remark}

%\begin{remark}
%We use the absolute value of $\mu_{n}^{\star}$ above because, in principle, there is no guarantee that $\mu_{n}^{\star}:= 2\mu_{n}-\mu_{n-1}$ is positive everywhere. However, this will seldom be needed in practice, since a continuous-in-time flow evolution will usually result in $2\mu_{n} >\mu_{n-1}$, unless the time-step size  is very large.
%\end{remark}

\subsection{A first-order fractional-step scheme}
Our IMEX treatment of the viscous term is especially advantageous when also segregating the pressure since, in that case, the resulting scheme will only require the solution of \textsl{scalar} subproblems: one for the pressure and one for each velocity component. For that, we propose an IMEX modification of the incremental pressure-correction method by \citet{Deteix2022}, who originally considered a fully implicit treatment of the convective and viscous terms. While classical projection methods for variable-density flows \cite{Guermond2000} require solving two pressure Poisson equations per time step, this recent variant requires only one, which is why it is our framework of choice here.

The type of boundary condition plays a key role when designing fractional-step schemes. For the momentum equation \eqref{momentum2}, we will consider, for simplicity, a pure Dirichlet boundary condition
\begin{align*}
    \ve{u} = \mathbf{v} \ \ \text{on} \ \partial\Omega\times(0,T]&\, ,
\end{align*}
%\begin{align}
 %   \ve{u} = \mathbf{v} \ \ \text{on} \ \Gamma_D\times(0,T]&\, ,\\
  %  \ve{n}\times[(\mu\nabla\ve{u})\ve{n}]  = \ve{0} \ \ \text{and} \ \ \ve{u}\cdot\ve{n} = 0\ \ \text{on} \ \Gamma_S\times(0,T]&\, ,
%\end{align}
where $\mathbf{v}$ is known. For this setup, we propose the following method: 
%\begin{align}
    %\ve{u} &= \mathbf{v} \ \ \text{on} \ \Gamma_D\times(0,T]\, ,\\
    %(\mu\nabla\ve{u})\ve{n} - p\ve{n} &= \ve{0} \ \ \text{on} \ \Gamma_N\times(0,T]\, ,
%\end{align}
%where $\Gamma_N$ and $\Gamma_D\neq \emptyset$ form a non-overlapping partition of $\partial\Omega$, and $\mathbf{v}$ is known. For this setup, we propose the following method:
\begin{itemize}
\item \textbf{Step 0:} To initialise the algorithm, set $\hat{\mathbf{u}}_0 = \ve{u}_0$ and provide (an approximation for) $p_0$.

\item \textbf{Step 1:} Find $\rho_{n+1}$ as the solution of 
\begin{equation}
\frac{\rho_{n+1}-\rho_n}{\tau} + \hat{\ve{u}}_n\cdot\nabla\rho_{n+1} + \frac{\nabla\cdot\hat{\ve{u}}_{n}}{2}\rho_{n+1} = 0\, ,
\end{equation}
with suitable boundary conditions, and update the viscosity $\mu$.

\item \textbf{Step 2:} Find the velocity $\ve{u}_{n+1}$ by solving
\begin{flalign}
\begin{cases}
\frac{\rho_{n+1}}{\tau}\ve{u}_{n+1} + (\rho_{n+1}\nabla\ve{u}_{n+1})\ve{u}_n + \frac{1}{2}[\nabla\cdot(\rho_{n+1}\ve{u}_n)]\ve{u}_{n+1}  -\nabla\cdot\left(\mu_{n+1}\nabla\ve{u}_{n+1}\right)  \\
 \hspace{16 mm} 
  = \frac{\sqrt{\rho_{n+1}\rho_n}}{\tau}\hat{\mathbf{u}}_{n} - \sqrt{\frac{\rho_{n+1}}{\rho_n}}\nabla p_{n} + \nabla^{\top}\ve{u}_{n}\nabla\mu_{n+1} +  \ve{f}_{n+1}\, ,\\
\ve{u}_{n+1}|_{\partial\Omega} = \mathbf{v}_{n+1}\, .
\end{cases} && \label{FS-momentum}
\end{flalign}

\item \textbf{Step 3:} Update the pressure through the Neumann problem
\begin{flalign}
\begin{cases}
-\nabla\cdot\big(\frac{1}{\rho_{n+1}}\nabla p_{n+1}\big) =  -\nabla\cdot\big(\frac{1}{\sqrt{\rho_n\rho_{n+1}}}\nabla p_{n}\big) - \frac{1}{\tau}\nabla\cdot\ve{u}_{n+1}\, ,\\
\left.\left[\ve{n}\cdot\Big(\frac{1}{\sqrt{\rho_{n+1}}}\nabla p_{n+1} - \frac{1}{\sqrt{\rho_{n}}}\nabla p_{n}\Big)\right]\right|_{\partial\Omega} = 0 \, .
\end{cases}&& \label{Poisson}
\end{flalign}

\item \textbf{Step 4:} Update the end-of-step velocity $\hat{\mathbf{u}}$ via 
\begin{flalign}
\hat{\mathbf{u}}_{n+1} &= \ve{u}_{n+1} - \frac{\tau}{\sqrt{\rho_{n+1}}}\left(\frac{1}{\sqrt{\rho_{n+1}}}\nabla p_{n+1}-\frac{1}{\sqrt{\rho_{n}}}\nabla p_{n}\right) \, . && \label{projectedU}
\end{flalign}
\end{itemize}

Notice that combining Eqs.~\eqref{Poisson} and \eqref{projectedU} implies, as in standard projection methods,
\begin{align}
\nabla\cdot\hat{\mathbf{u}}_{n+1} = 0 \ \ \text{and}
    \ \ \ve{n}\cdot(\ve{u}_{n+1}-\hat{\mathbf{u}}_{n+1})|_{\partial\Omega} = 0\, .\label{divEndOfStep}
\end{align}

It is simple to verify that Step 2 can be solved individually for each $i$th velocity component $u_{n+1}^{i}$:
\begin{equation}
    \begin{split}
       \frac{\rho_{n+1}}{\tau}u^i_{n+1} + \rho_{n+1}\ve{u}_n\cdot\nabla u^i_{n+1} + \frac{1}{2}[\nabla\cdot(\rho_{n+1}\ve{u}_n)]
u^i_{n+1}  -\nabla\cdot\left(\mu_{n+1}\nabla u^i_{n+1}\right)  \\
 \hspace{16 mm} 
  = \frac{\sqrt{\rho_{n+1}\rho_n}}{\tau}\hat{\mathrm{u}}^i_{n} - \sqrt{\frac{\rho_{n+1}}{\rho_n}}\frac{\partial p_{n}}{\partial x_i} + \frac{\partial \ve{u}_n}{\partial x_i}\cdot\nabla\mu_{n+1} +  f^i_{n+1}\, , \quad i=1,...,d\, .
  \label{componentsGU} 
    \end{split}
\end{equation}
Our fractional-step algorithm is thus very simple to implement and only requires solving (in addition to the density equation) $d$ linear scalar transport equations plus one Poisson-like problem. Moreover, up to boundary conditions, the system matrices for each of the $i$th velocity equations are identical, which reduces assembling efforts (see Section \ref{sec_impl}).

\begin{remark}
    Since $\ve{n}\cdot\hat{\mathbf{u}}_{n+1} = \ve{n}\cdot\ve{u}_{n+1}$ on $\partial\Omega$, we can use either $\hat{\mathbf{u}}_{n+1}$ or $\ve{u}_{n+1}$ as convective velocities for the density and momentum equations, without losing the skew-symmetry properties \eqref{skewSymRho} and \eqref{skewSymU}. However, enforcing Eq.~\eqref{projectedU} strongly results in $\hat{\mathbf{u}}_{n+1}$ less regular than $\ve{u}_{n+1}$, so computing derivatives of $\hat{\mathbf{u}}_{n+1}$ in the fully discrete case may require some regularisation. When using finite elements, for example, a simple $L^2(\Omega)$ projection of $\hat{\mathbf{u}}_{n+1}$ onto a continuous finite element space will suffice. 
\end{remark}
\begin{remark}
    To extend the projection step to the more general case where both Neumann and Dirichlet conditions are considered for the momentum equation, the Neumann boundary condition in \eqref{Poisson} must be partially replaced by a Dirichlet one (see, e.g., Section 10 in the seminal overview by \citet{Guermond2006}).
\end{remark}

\begin{remark}
The stability analysis of the scheme just presented will be performed independently of the approximation $p_0$. In fact, the stability analysis will be completely independent of the concrete approximation. A proposal for the approximation of $p(0)$ guaranteeing optimal error estimates was presented by \citet{JN15}.
\end{remark}

\subsection{On second-order fractional-step schemes}
In principle, the fractional-step algorithm presented above could be extended to second order by combining, for instance, BDF2 and second-order extrapolations. However, as remarked by \citet{Wu2017} -- and still true, to the best of our knowledge --, there is still no provably stable second-order projection method for the variable-density problem. In this context, we shall presently refrain from proposing an IMEX fractional-step scheme of second order. 

\section{Temporal stability analysis}\label{sec_analysis}
We will next derive temporal stability proofs for each of the three schemes proposed herein (a reader not so interested in the theory may skip to Section \ref{sec_impl}, where we discuss implementation). As usual for the analysis, we assume homogeneous Dirichlet data: $\ve{u}_{n+1} = \ve{0}$ on $\partial\Omega$ for all $t_{n+1}$. In the proofs, $C$ denotes a generic positive constant depending only on initial conditions, on problem data and on $\tau$, and which does not grow with decreasing $\tau$ ($\partial_{\tau} C \geq 0$ for all $\tau > 0$). Although we will not discuss spatial stability, the theory presented in this section holds, e.g., for standard $H^1$-conforming finite element spaces. Before analysing the velocity, we state the stability of $\rho$.
\begin{lemma}[Density stability] For any time-step size $\tau= T/N$, the BDF schemes \eqref{densityBDF1} and \eqref{densityBDF2} yield, respectively,
\begin{align*}
\|\rho_{N}\|^2 + \sum_{n=0}^{N}\|\delta\rho_{n}\|^2 = \|\rho_{0}\|^2
\end{align*}
and
\begin{align*}
\|\rho_{N}\|^2 + \|\rho_{N}^{\star}\|^2 + \sum_{n=1}^{N}\|\delta^2\rho_{n}\|^2 = \|\rho_{1}\|^2 + \|\rho_{1}^{\star}\|^2\, .
\end{align*} 
\end{lemma}
\proof{These classical results follow immediately from the skew-symmetry \eqref{skewSymRho} and the identities \eqref{identityBDF1} and \eqref{identityBDF2}.
}

%Here, we consider Eqs.~\eqref{mass} or \eqref{mass2} as a prototypical way to propagate the density field. Other possibilities can be used instead, since the velocity/pressure stability proofs in the next section are derived without invoking the density equation. In particular, for a two-phase flow, we can replace the density problem by the level-set equation
%\begin{equation*}
    %\partial_t\phi + \ve{u}\cdot\nabla\phi + %\frac{\nabla\cdot\ve{u}}{2}\phi = 0\, ,
%\end{equation*}
%where $\phi$ is a signed distance function from which we update $\rho$ and $\nu$ (see Section \ref{sec_damBreak}). Other examples are phase-field equations such as the Cahn--Hilliard and Allen--Cahn models \cite{Shen2010}.

\begin{remark}
    The following stability estimates assume that at each time step $t=t_{n+1}$, the forcing term satisfies $\ve{f}_{n+1}\in [H^{-1}(\Omega)]^d$. Modified estimates considering $\ve{f}\in [L^2(Q)]^d$ can also be derived, as done in our recent work \cite{Barrenechea2024}.
\end{remark}

\subsection{First-order coupled scheme}
For the first-order coupled IMEX scheme, we will prove the following stability estimate. 
\begin{theorem}[Stability of the first-order coupled scheme]
For any time-step size $\tau=T/N$, $N\ge 1$, the IMEX scheme (\eqref{densityBDF1},\eqref{coupledBDF1},\eqref{divuBDF1}) satisfies the stability estimate
    \begin{equation}
    \begin{split}      
    &\|\sqrt{\rho_N}\,\ve{u}_{N}\|^2 + \tau\varepsilon\|\sqrt{\mu_N}\,\nabla\ve{u}_{N}\|^2  + \tau(2-\epsilon-\varepsilon)\sum_{n=1}^{N}\big\|\sqrt{\mu_{n}}\,\nabla\ve{u}_{n}\big\|^2 \\ 
     &\leq \left(C + \sum_{n=1}^{N}\frac{\tau}{\epsilon\mu_{\mathrm{min}}}\|\ve{f}_{n}\|_{-1}^2\right)\mathrm{exp}\left(\frac{ T\|\nabla\mu \|_{\infty}^2}{\varepsilon\varrho_{\mathrm{min}}\mu_{\mathrm{min}}}\right)
     \, ,\label{stabilityBDF1}
     \end{split}
    \end{equation}
where $\varepsilon,\epsilon$ are two positive constants such that $\varepsilon + \epsilon \leq 2$.
\end{theorem}

\proof{We start by testing Eqs.~\eqref{coupledBDF1} and \eqref{divuBDF1} with $2\tau\ve{u}_{n+1}$ and $2\tau p_{n+1}$, respectively, and adding the two results, which yields
\begin{equation}
     \|\ve{U}_{n+1}\|^2   - \|\ve{U}_{n}\|^2+ \|\delta\ve{U}_{n+1}\|^2 + 2\tau\big\|\sqrt{\mu_{n+1}}\,\nabla\ve{u}_{n+1}\big\|^2 = 2\tau(\nabla^{\top}\ve{u}_{n}\nabla\mu_{n+1},\ve{u}_{n+1}) + 2\tau\langle\ve{f}_{n+1},\ve{u}_{n+1}\rangle ,\label{CC}
\end{equation}
due to \eqref{identityBDF1} and \eqref{divuBDF1}.
We must bound both terms on the right-hand side of \eqref{CC}. Using Hölder's and Young's inequalities, we get 
\begin{align}
2\tau(\nabla^{\top}\ve{u}_{n}\nabla\mu_{n+1},\ve{u}_{n+1}) &= 2\tau \bigg(\frac{1}{\sqrt{\rho_{n+1}}}\nabla^{\top}\ve{u}_{n}\nabla\mu_{n+1},\ve{U}_{n+1}\bigg)\nonumber \\
&=2\tau \bigg(\frac{1}{\sqrt{\rho_{n+1}}}\nabla^{\top}\ve{u}_{n}\nabla\mu_{n+1},\delta\ve{U}_{n+1} + \ve{U}_n\bigg) \nonumber\\
 &\leq 2\tau \frac{\|\nabla\mu\|_{\infty}}{\sqrt{\varrho_{\mathrm{min}}}} \|\nabla\ve{u}_n \|\,\|\delta\ve{U}_{n+1} \|  + 2\tau\frac{\|\nabla\mu\|_{\infty}}{\sqrt{\varrho_{\mathrm{min}}}} \|\nabla\ve{u}_n \|\, \|\ve{U}_{n} \| \nonumber\\
 &\leq \|\delta\ve{U}_{n+1} \|^2 + \tau^2\frac{\|\nabla\mu\|_{\infty}^2}{\varrho_{\mathrm{min}}} \|\nabla\ve{u}_n \|^2 + \tau\varepsilon\mu_{\mathrm{min}} \|\nabla\ve{u}_n \|^2  + \frac{\tau\|\nabla\mu\|_{\infty}^2}{\varepsilon\varrho_{\mathrm{min}}\mu_{\mathrm{min}}} \|\ve{U}_{n}\|^2 \nonumber \\
 &= \|\delta\ve{U}_{n+1} \|^2 +\tau\varepsilon\mu_{\mathrm{min}} \|\nabla\ve{u}_n \|^2 +\frac{\tau\|\nabla\mu\|_{\infty}^2}{\varepsilon\varrho_{\mathrm{min}}\mu_{\mathrm{min}}} (\|\ve{U}_{n}\|^2 + \tau\varepsilon\mu_{\mathrm{min}}\|\nabla\ve{u}_{n}\|^2)\nonumber\\
 &\leq \|\delta\ve{U}_{n+1} \|^2 +\tau\varepsilon \|\sqrt{\mu_n}\,\nabla\ve{u}_n \|^2 +\frac{\tau\|\nabla\mu\|_{\infty}^2}{\varepsilon\varrho_{\mathrm{min}}\mu_{\mathrm{min}}} (\|\ve{U}_{n}\|^2 + \tau\varepsilon\|\sqrt{\mu_n}\,\nabla\ve{u}_{n}\|^2)
 \, ,  \label{AA}
\end{align}
and
\begin{align}
    2\tau\langle\ve{f}_{n+1},\ve{u}_{n+1}\rangle 
    &\leq 2\tau\|\ve{f}_{n+1}\|_{-1}\|\nabla\ve{u}_{n+1}\| \nonumber\\ 
   &\leq \frac{\tau}{\epsilon\mu_{\mathrm{min}}}\|\ve{f}_{n+1}\|_{-1}^2 + \tau\epsilon\mu_{\mathrm{min}}\|\nabla\ve{u}_{n+1}\|^2 \nonumber \\
   &\leq \frac{\tau}{\epsilon\mu_{\mathrm{min}}}\|\ve{f}_{n+1}\|_{-1}^2 + \tau\epsilon\|\sqrt{\mu_{n+1}}\,\nabla\ve{u}_{n+1}\|^2\,,\label{BB}
\end{align}
for arbitrary $\varepsilon,\epsilon>0$. Inserting \eqref{AA} and \eqref{BB} in \eqref{CC} gives
\begin{align*}
     &\|\ve{U}_{n+1}\|^2 - \|\ve{U}_{n}\|^2 + \tau(2-\epsilon)\big\|\sqrt{\mu_{n+1}}\,\nabla\ve{u}_{n+1}\big\|^2 \\ 
     &\leq \tau\varepsilon\|\sqrt{\mu_n}\,\nabla\ve{u}_{n}\|^2+ \frac{\tau\|\nabla\mu\|_{\infty}^2}{\varepsilon\varrho_{\mathrm{min}}\mu_{\mathrm{min}}} (\|\ve{U}_{n}\|^2 + \tau\varepsilon\|\sqrt{\mu_n}\,\nabla\ve{u}_{n}\|^2) + \frac{\tau}{\epsilon\mu_{\mathrm{min}}}\|\ve{f}_{n+1}\|_{-1}^2
     \, .
\end{align*}
Adding up from $n=0$ to $n=N-1$ yields
\begin{align*}
     &\|\ve{U}_{N}\|^2 + \tau\varepsilon\|\sqrt{\mu_N}\,\nabla\ve{u}_{N}\|^2  + \tau(2-\epsilon-\varepsilon)\sum_{n=1}^{N}\big\|\sqrt{\mu_{n}}\,\nabla\ve{u}_{n}\big\|^2 \\ 
     &\leq \|\ve{U}_{0}\|^2 +\tau\varepsilon\|\sqrt{\mu_0}\,\nabla\ve{u}_{0}\|^2+ \frac{\tau\|\nabla\mu\|_{\infty}^2}{\varepsilon\varrho_{\mathrm{min}}\mu_{\mathrm{min}}} \sum_{n=0}^{N-1}(\|\ve{U}_{n}\|^2 + \tau\varepsilon\|\sqrt{\mu_n}\,\nabla\ve{u}_{n}\|^2) + \sum_{n=1}^{N}\frac{\tau}{\epsilon\mu_{\mathrm{min}}}\|\ve{f}_{n}\|_{-1}^2
     \, .
\end{align*}
Estimate \eqref{stabilityBDF1} now follows directly from the discrete Gronwall lemma, with $a_n = \|\ve{U}_{n}\|^2 + \tau\varepsilon\|\sqrt{\mu_n}\,\nabla\ve{u}_{n}\|^2$, $b_n = \tau(2-\epsilon-\varepsilon)\big\|\sqrt{\mu_{n}}\,\nabla\ve{u}_{n}\big\|^2$ and $c_n=\frac{\tau}{\epsilon\mu_{\mathrm{min}}}\|\ve{f}_{n}\|_{-1}^2$.}

\subsection{Second-order coupled scheme}
Unlike the scheme just analysed, the second-order version is not self-starting, since $\rho_{n-1}$ and $\ve{u}_{n-1}$ are not defined at the first time step ($n=0$). The simplest approach is to start with the first-order method (\eqref{densityBDF1},\eqref{coupledBDF1},\eqref{divuBDF1}). From the second time step onwards we use the second-order scheme, for which we will prove the following stability result.
\begin{theorem}[Stability of the second-order coupled scheme]
For any time-step size $\tau=T/N$, $N\ge 2$, the IMEX scheme (\eqref{densityBDF2},\eqref{coupledBDF2},\eqref{divuBDF2}) satisfies the stability estimate
    \begin{equation}
    \begin{split}      
    &\big\|\sqrt{\rho_N}\,\ve{u}_{N}\big\|^2  + \big\|2\sqrt{\rho_N}\,\ve{u}_{N}-\sqrt{\rho_{N-1}}\,\ve{u}_{N-1}\big\|^2 + 40\tau\varepsilon\big\|\sqrt{\mu_N}\,\nabla\ve{u}_{N}\big\|^2  \\
    &+ 4\tau\sum_{n=1}^{N}\left[(1-\epsilon-10\varepsilon)\big\|\sqrt{\mu_{n}}\,\nabla\ve{u}_{n}\big\|^2 + \varepsilon\mu_{\mathrm{min}}\big\|2\nabla\ve{u}_{n-1}+\nabla\ve{u}_{n-2}\big\|^2\right]  \\   
     &\leq \left(\tilde{C} + \sum_{n=1}^{N}\frac{\tau}{\epsilon\mu_{\mathrm{min}}}\|\ve{f}_{n}\|_{-1}^2\right)\mathrm{exp}\left(\frac{ T\|\nabla\mu \|_{\infty}^2}{\varepsilon\varrho_{\mathrm{min}}\mu_{\mathrm{min}}}\right)
     \, ,\label{stabilityBDF2}
     \end{split}
    \end{equation}
where $\varepsilon,\epsilon$ are two positive constants such that $10\varepsilon + \epsilon \leq 1$.
\end{theorem}

\proof{We start by testing Eqs.~\eqref{coupledBDF2} and \eqref{divuBDF2} with $4\tau\ve{u}_{n+1}$ and $4\tau p_{n+1}$, respectively, and adding those results, so that
\begin{equation}
\begin{split}
     \|\ve{U}_{n+1}\|^2 + \|\ve{U}_{n+1}^{\star}\|^2   - \|\ve{U}_{n}\|^2 - \|\ve{U}_{n}^{\star}\|^2 + \|\delta^2\ve{U}_{n+1}\|^2 + 4\tau\big\|\sqrt{\mu_{n+1}}\,\nabla\ve{u}_{n+1}\big\|^2 \\
     = 4\tau(\nabla^{\top}\ve{u}_{n}^{\star}\nabla\mu_{n+1},\ve{u}_{n+1}) + 4\tau\langle\ve{f}_{n+1},\ve{u}_{n+1}\rangle\,.\label{CCC}
     \end{split}
\end{equation} 
The last term on the right-hand side is estimated as:
\begin{align}
    4\tau\langle\ve{f}_{n+1},\ve{u}_{n+1}\rangle 
   \leq \frac{\tau}{\epsilon\mu_{\mathrm{min}}}\|\ve{f}_{n+1}\|_{-1}^2 + 4\tau\epsilon\|\sqrt{\mu_{n+1}}\,\nabla\ve{u}_{n+1}\|^2\,.\label{BBB}
\end{align}
For the remaining term we have
\begin{align*}
4\tau(\nabla^{\top}\ve{u}_{n}^{\star}\nabla\mu_{n+1},\ve{u}_{n+1}) &= 4\tau \bigg(\frac{1}{\sqrt{\rho_{n+1}}}\nabla^{\top}\ve{u}_{n}^{\star}\nabla\mu_{n+1},\ve{U}_{n+1}\bigg) \\
&=4\tau \bigg(\frac{1}{\sqrt{\rho_{n+1}}}\nabla^{\top}\ve{u}_{n}^{\star}\nabla\mu_{n+1},\delta^2\ve{U}_{n+1} + \ve{U}_n^{\star}\bigg) \\
 &\leq 4\tau \frac{\|\nabla\mu\|_{\infty}}{\sqrt{\varrho_{\mathrm{min}}}} \|\nabla\ve{u}_n^{\star} \|\,\|\delta^2\ve{U}_{n+1} \|  + 4\tau\frac{\|\nabla\mu\|_{\infty}}{\sqrt{\varrho_{\mathrm{min}}}} \|\nabla\ve{u}_n^{\star} \|\, \|\ve{U}_{n}^{\star} \| \\
 &\leq \|\delta^2\ve{U}_{n+1} \|^2 + 4\tau^2\frac{\|\nabla\mu\|_{\infty}^2}{\varrho_{\mathrm{min}}} \|\nabla\ve{u}_n^{\star} \|^2 + 4\tau\varepsilon\mu_{\mathrm{min}} \|\nabla\ve{u}_n^{\star} \|^2  + \frac{\tau\|\nabla\mu\|_{\infty}^2}{\varepsilon\varrho_{\mathrm{min}}\mu_{\mathrm{min}}} \|\ve{U}_{n}^{\star}\|^2 .
\end{align*}
For the viscous terms above, we use
\begin{align*}
\mu_{\text{min}}\|\nabla\ve{u}_n^{\star}\|^2 &= \mu_{\text{min}}\|2\nabla\ve{u}_n - \nabla\ve{u}_{n-1}\|^2\\
&\equiv \mu_{\text{min}}(2\|2\nabla\ve{u}_n\|^2 + 2\|\nabla\ve{u}_{n-1}\|^2-\|2\nabla\ve{u}_n + \nabla\ve{u}_{n-1}\|^2)\\
    &\leq 2( \|2\sqrt{\mu_{n}}\,\nabla\ve{u}_n\|^2 + \|\sqrt{\mu_{n-1}}\,\nabla\ve{u}_{n-1}\|^2)-\mu_{\text{min}}\|2\nabla\ve{u}_n + \nabla\ve{u}_{n-1}\|^2\\
    & \leq 8\big\|\sqrt{\mu_{n}}\,\nabla\ve{u}_n\big\|^2 + 2\big\|\sqrt{\mu_{n-1}}\,\nabla\ve{u}_{n-1}\big\|^2\, ,
\end{align*}
so that
\begin{equation}
    \begin{split}     &4\tau(\nabla^{\top}\ve{u}_{n}^{\star}\nabla\mu_{n+1},\ve{u}_{n+1}) \\
    &\leq \|\delta^2\ve{U}_{n+1} \|^2 + \tau\varepsilon\left(32\big\|\sqrt{\mu_{n}}\,\nabla\ve{u}_n\big\|^2 + 8\big\|\sqrt{\mu_{n-1}}\,\nabla\ve{u}_{n-1}\big\|^2 - 4\mu_{\text{min}}\big\|2\nabla\ve{u}_{n}+\nabla\ve{u}_{n-1}\big\|^2\right)\\
    &+ \frac{\tau\|\nabla\mu\|_{\infty}^2}{\varepsilon\varrho_{\mathrm{min}}\mu_{\mathrm{min}}}\left(\|\ve{U}_{n}^{\star} \|^2 + 32\varepsilon\tau\big\|\sqrt{\mu_{n}}\,\nabla\ve{u}_n\big\|^2 + 8\varepsilon\tau\big\|\sqrt{\mu_{n-1}}\,\nabla\ve{u}_{n-1}\big\|^2\right).
    \end{split}
    \label{AAA}
\end{equation}
Combining Eqs.~\eqref{CCC}, \eqref{BBB} and \eqref{AAA} yields
\begin{equation}
\begin{split}
     &\|\ve{U}_{n+1}\|^2 + \|\ve{U}_{n+1}^{\star}\|^2   + 4\tau(1-\epsilon)\big\|\sqrt{\mu_{n+1}}\,\nabla\ve{u}_{n+1}\big\|^2 + 4\tau\varepsilon\mu_{\text{min}}\big\|2\nabla\ve{u}_{n}+\nabla\ve{u}_{n-1}\big\|^2\\
     &\leq \|\ve{U}_{n}\|^2 + \|\ve{U}_{n}^{\star}\|^2 + \frac{\tau}{\epsilon\mu_{\mathrm{min}}}\|\ve{f}_{n+1}\|_{-1}^2 +  4\tau\varepsilon\left(8\|\sqrt{\mu_{n}}\,\nabla\ve{u}_n\|^2 + 2\|\sqrt{\mu_{n-1}}\,\nabla\ve{u}_{n-1}\|^2\right) \\
     &+ \frac{\tau\|\nabla\mu\|_{\infty}^2}{\varepsilon\varrho_{\mathrm{min}}\mu_{\mathrm{min}}}\left(\|\ve{U}_{n}^{\star} \|^2 + 32\varepsilon\tau\|\sqrt{\mu_{n}}\,\nabla\ve{u}_n\|^2 + 8\varepsilon\tau\|\sqrt{\mu_{n-1}}\,\nabla\ve{u}_{n-1}\|^2\right)\,.\label{DDD}
     \end{split}
\end{equation} 

We will now sum \eqref{DDD} from $n=1$ to $N-1$. For the viscous term, we have that
\begin{align*}
        &\sum_{n=1}^{N-1}\left[(1-\epsilon)\big\|\sqrt{\mu_{n+1}}\,\nabla\ve{u}_{n+1}\big\|^2 -  8\varepsilon\|\sqrt{\mu_{n}}\,\nabla\ve{u}_n\|^2 -2\varepsilon\|\sqrt{\mu_{n-1}}\,\nabla\ve{u}_{n-1}\|^2\right]\\
        &= 10\varepsilon\left(\big\|\sqrt{\mu_{N}}\,\nabla\ve{u}_{N}\big\|^2 - \|\sqrt{\mu_{1}}\,\nabla\ve{u}_{1}\big\|^2\right) + 2\varepsilon\left(\big\|\sqrt{\mu_{N-1}}\,\nabla\ve{u}_{N-1}\big\|^2 - \|\sqrt{\mu_{0}}\,\nabla\ve{u}_{0}\big\|^2\right) \\
        &+
        (1-\epsilon-10\varepsilon)\sum_{n=1}^{N-1}\big\|\sqrt{\mu_{n+1}}\,\nabla\ve{u}_{n+1}\big\|^2\\
        &\geq 10\varepsilon\big\|\sqrt{\mu_{N}}\,\nabla\ve{u}_{N}\big\|^2-\varepsilon\Big(10\|\sqrt{\mu_{1}}\,\nabla\ve{u}_{1}\big\|^2 + 2\|\sqrt{\mu_{0}}\,\nabla\ve{u}_{0}\big\|^2 \Big)+
        (1-\epsilon-10\varepsilon)\sum_{n=1}^{N-1}\big\|\sqrt{\mu_{n+1}}\,\nabla\ve{u}_{n+1}\big\|^2,
\end{align*}
hence
\begin{equation}
\begin{split}
     &\|\ve{U}_{N}\|^2 + \|\ve{U}_{N}^{\star}\|^2  + 40\varepsilon\tau\big\|\sqrt{\mu_{N}}\,\nabla\ve{u}_{N}\big\|^2 \\
        &+4\tau\sum_{n=2}^{N}\left[(1-\epsilon-10\varepsilon)\big\|\sqrt{\mu_{n}}\,\nabla\ve{u}_{n}\big\|^2 + \varepsilon\mu_{\mathrm{min}}\big\|2\nabla\ve{u}_{n-1}+\nabla\ve{u}_{n-2}\big\|^2\right] \\
     &\leq  \|\ve{U}_{1}\|^2 + \|\ve{U}_{1}^{\star}\|^2 + 4\tau\varepsilon\left(10\|\sqrt{\mu_{1}}\,\nabla\ve{u}_{1}\big\|^2 + 2\|\sqrt{\mu_{0}}\,\nabla\ve{u}_{0}\big\|^2 \right) + \frac{\tau}{\epsilon\mu_{\mathrm{min}}}\sum_{n=2}^{N}\|\ve{f}_{n}\|_{-1}^2  \\
     &+ \frac{\tau\|\nabla\mu\|_{\infty}^2}{\varepsilon\varrho_{\mathrm{min}}\mu_{\mathrm{min}}}\sum_{n=1}^{N-1}\left(\|\ve{U}_{n}^{\star} \|^2 + 32\varepsilon\tau\|\sqrt{\mu_{n}}\,\nabla\ve{u}_n\|^2 + 8\varepsilon\tau\|\sqrt{\mu_{n-1}}\,\nabla\ve{u}_{n-1}\|^2\right)\\
     &\leq C + \frac{\tau\|\nabla\mu\|_{\infty}^2}{\varepsilon\varrho_{\mathrm{min}}\mu_{\mathrm{min}}}\sum_{n=1}^{N-1}\left(\|\ve{U}_{n}^{\star} \|^2 + 40\varepsilon\tau\|\sqrt{\mu_{n}}\,\nabla\ve{u}_n\|^2\right) +  \frac{\tau}{\epsilon\mu_{\mathrm{min}}}\sum_{n=2}^{N}\|\ve{f}_{n}\|_{-1}^2\\
     &\leq \tilde{C} + \frac{\tau\|\nabla\mu\|_{\infty}^2}{\varepsilon\varrho_{\mathrm{min}}\mu_{\mathrm{min}}}\sum_{n=1}^{N-1}\left(\|\ve{U}_{n}\|^2 + \|\ve{U}_{n}^{\star} \|^2 + 40\varepsilon\tau\|\sqrt{\mu_{n}}\,\nabla\ve{u}_n\|^2\right) +  \frac{\tau}{\epsilon\mu_{\mathrm{min}}}\sum_{n=1}^{N}\|\ve{f}_{n}\|_{-1}^2
     \,,\label{EEE}
     \end{split}
\end{equation}
from which stability follows provided that $\epsilon + 10\varepsilon \leq 1$, using the Gronwall inequality.
}
\begin{remark}
In \eqref{EEE}, we were able to incorporate the contributions from the first time step into the finite constant $\tilde{C}$ because of the already proven stability of the first-order scheme, which is used for initialisation. 
\end{remark}

\subsection{Fractional-step scheme}
We will now prove that our fractional-step IMEX method is also unconditionally stable.
\begin{theorem}[Stability of the fractional-step IMEX scheme]
For any time-step size $\tau=T/N$, $N\ge 1$, the IMEX scheme \eqref{FS-momentum}--\eqref{projectedU} satisfies the stability estimate
    \begin{equation}
    \begin{split}      
    &\|\sqrt{\rho_N}\,\hat{\mathbf{u}}_{N}\|^2 + \tau\varepsilon\|\sqrt{\mu_N}\,\nabla\ve{u}_{N}\|^2  + \tau^2\left\|\frac{1}{\sqrt{\rho_{N}}}\nabla p_{N}\right\|^2 + (2-\epsilon-\varepsilon)\tau\sum_{n=1}^{N}\big\|\sqrt{\mu_{n}}\,\nabla\ve{u}_{n}\big\|^2  \\ 
     &\leq \Big(C + \sum_{n=1}^{N}\frac{\tau}{\epsilon\mu_{\mathrm{min}}}\|\ve{f}_{n}\|_{-1}^2\Big)\mathrm{exp}\Big(\frac{ T\|\nabla\mu \|_{\infty}^2}{\varepsilon\varrho_{\mathrm{min}}\mu_{\mathrm{min}}}\Big)
     \, ,\label{stabilityGU}
     \end{split}
    \end{equation}
where $\varepsilon,\epsilon$ are two positive constants such that $\varepsilon + \epsilon \leq 2$.
\end{theorem}
\proof{Let us test the first equation in \eqref{FS-momentum} against $2\tau\ve{u}_{n+1}$, so that
\begin{equation}
\begin{split}
     &\|\ve{U}_{n+1}\|^2   - \|\hat{\mathbf{U}}_{n}\|^2+ \|\ve{U}_{n+1} - \hat{\mathbf{U}}_{n}\|^2 + 2\tau\big\|\sqrt{\mu_{n+1}}\,\nabla\ve{u}_{n+1}\big\|^2 \\
    & = 2\tau(\nabla^{\top}\ve{u}_{n}\nabla\mu_{n+1},\ve{u}_{n+1}) + 2\tau\langle\ve{f}_{n+1},\ve{u}_{n+1}\rangle - 2\tau\left(\chi_n\nabla p_n,\ve{U}_{n+1}\right) ,\label{proofFS1}
     \end{split}
\end{equation}
introducing the short-hand notation $\chi_n:= 1/\sqrt{\rho_n}$. Then, testing \eqref{Poisson} with $2\tau p_{n+1}$ gives
\begin{align*}
        \tau\big(\|\chi_{n+1}\nabla p_{n+1}\|^2 -  \|\chi_{n}\nabla p_{n}\|^2 + \|\chi_{n+1}\nabla p_{n+1} - \chi_{n}\nabla p_{n}\|^2\big) &= 2 (-\nabla\cdot\ve{u}_{n+1},p_{n+1})\\
        &= 2 (\nabla\cdot\hat{\mathbf{u}}_{n+1}-\nabla\cdot\ve{u}_{n+1},p_{n+1})\\
        &= 2 (\ve{u}_{n+1}-\hat{\mathbf{u}}_{n+1},\nabla p_{n+1})\\
        &= 2 (\ve{U}_{n+1}-\hat{\mathbf{U}}_{n+1},\chi_{n+1}\nabla p_{n+1})\, ,
\end{align*}
where we have used \eqref{divEndOfStep}. This, combined with  Eq.~\eqref{projectedU}, yields
\begin{align}
        \tau^2\|\chi_{n+1}\nabla p_{n+1}\|^2 -  \tau^2\|\chi_{n}\nabla p_{n}\|^2 + \|\ve{U}_{n+1} - \hat{\mathbf{U}}_{n}\|^2 =  2\tau(\ve{U}_{n+1}-\hat{\mathbf{U}}_{n+1},\chi_{n+1}\nabla p_{n+1})\, . \label{proofFS2}
\end{align}

We now test Eq.~\eqref{projectedU} against $2\rho_{n+1}\hat{\mathbf{u}}_{n+1}$, resulting in
\begin{equation}
\begin{split}
    \|\hat{\mathbf{U}}_{n+1}\|^2 - \|\ve{U}_{n+1}\|^2   + \|\ve{U}_{n+1} - \hat{\mathbf{U}}_{n+1}\|^2 &= -2\tau (\chi_{n+1}\nabla p_{n+1} - \chi_{n}\nabla p_{n},\hat{\mathbf{U}}_{n+1})\\
    &= 2\tau (\chi_{n}\nabla p_{n},\hat{\mathbf{U}}_{n+1}) - 2\tau (\nabla p_{n+1},\hat{\mathbf{u}}_{n+1})\\
    &= 2\tau (\chi_{n}\nabla p_{n},\hat{\mathbf{U}}_{n+1}) \, , \label{proofFS3}
\end{split}
\end{equation}
again due to \eqref{divEndOfStep}. So, adding Eqs.\eqref{proofFS1}--\eqref{proofFS3} yields 
\begin{equation}
\begin{split}
     &\|\hat{\mathbf{U}}_{n+1}\|^2   - \|\hat{\mathbf{U}}_{n}\|^2+ \|\ve{U}_{n+1} - \hat{\mathbf{U}}_{n}\|^2 + 2\tau\big\|\sqrt{\mu_{n+1}}\,\nabla\ve{u}_{n+1}\big\|^2 + \tau^2\big(\|\chi_{n+1}\nabla p_{n+1}\|^2 -  \|\chi_{n}\nabla p_{n}\|^2\big) \\
    &= 2\tau(\nabla^{\top}\ve{u}_{n}\nabla\mu_n,\ve{u}_{n+1}) + 2\tau\langle\ve{f}_{n+1},\ve{u}_{n+1}\rangle \\
    & + 2\tau(\chi_{n+1}\nabla p_{n+1}-\chi_n\nabla p_n,\ve{U}_{n+1}- \hat{\mathbf{U}}_{n+1}) - 2\|\ve{U}_{n+1} - \hat{\mathbf{U}}_{n+1}\|^2 \, ,\label{proofFS4}
     \end{split}
\end{equation}
where the last two terms on the right-hand side cancel out, since
\begin{align*}
    \tau(\chi_{n+1}\nabla p_{n+1}-\chi_n\nabla p_n,\ve{U}_{n+1}- \hat{\mathbf{U}}_{n+1}) = \|\ve{U}_{n+1} - \hat{\mathbf{U}}_{n+1}\|^2 \, ,
\end{align*}
by construction (of the projected velocity $\hat{\mathbf{u}}_{n+1}$, see Eq.~\eqref{projectedU}).

For the explicit part of the viscous term in \eqref{proofFS4} we write, similarly as in \eqref{AA},
\begin{align}
2\tau(\nabla^{\top}\ve{u}_{n}\nabla\mu_{n+1},\ve{u}_{n+1}) &= 2\tau \left(\frac{1}{\sqrt{\rho_{n+1}}}\nabla^{\top}\ve{u}_{n}\nabla\mu_{n+1},\ve{U}_{n+1}\right)\nonumber \\
&=2\tau \left(\frac{1}{\sqrt{\rho_{n+1}}}\nabla^{\top}\ve{u}_{n}\nabla\mu_{n+1},\ve{U}_{n+1} - \hat{\mathbf{U}}_n + \hat{\mathbf{U}}_n\right) \nonumber\\
 &\leq 2\tau \frac{\|\nabla\mu\|_{\infty}}{\sqrt{\varrho_{\mathrm{min}}}} \|\nabla\ve{u}_n \|\,\|\ve{U}_{n+1} - \hat{\mathbf{U}}_{n} \|  + 2\tau\frac{\|\nabla\mu\|_{\infty}}{\sqrt{\varrho_{\mathrm{min}}}} \|\nabla\ve{u}_n \|\, \|\hat{\mathbf{U}}_{n} \| \nonumber\\
 &\leq \|\ve{U}_{n+1} - \hat{\mathbf{U}}_{n} \|^2 +\tau\varepsilon \|\sqrt{\mu_n}\,\nabla\ve{u}_n \|^2 +\frac{\tau\|\nabla\mu\|_{\infty}^2}{\varepsilon\varrho_{\mathrm{min}}\mu_{\mathrm{min}}} (\|\hat{\mathbf{U}}_{n}\|^2 + \tau\varepsilon\|\sqrt{\mu_n}\,\nabla\ve{u}_{n}\|^2)
 \, .  \label{explicitTermFS}
\end{align}

Finally, by estimating the forcing term as we did for the coupled scheme and combining that estimate with \eqref{proofFS4} and \eqref{explicitTermFS}, we obtain
\begin{align*}
     &\|\hat{\mathbf{U}}_{n+1}\|^2   - \|\hat{\mathbf{U}}_{n}\|^2  + \tau^2(\|\chi_{n+1}\nabla p_{n+1}\|^2 - \|\chi_{n}\nabla p_{n}\|^2) + (2-\epsilon)\tau\big\|\sqrt{\mu_{n+1}}\,\nabla\ve{u}_{n+1}\big\|^2  \\
    & \leq \varepsilon\tau\big\|\sqrt{\mu_{n}}\,\nabla\ve{u}_{n}\big\|^2 +  \frac{\tau\|\nabla\mu\|_{\infty}^2}{\varepsilon\varrho_{\mathrm{min}}\mu_{\mathrm{min}}} \left(\|\hat{\mathbf{U}}_{n}\|^2 + \tau\varepsilon\|\sqrt{\mu_n}\,\nabla\ve{u}_{n}\|^2\right)  + \frac{\tau}{\epsilon\mu_{\mathrm{min}}}\|\ve{f}_{n+1}\|_{-1}^2 \\
    &\leq \varepsilon\tau\big\|\sqrt{\mu_{n}}\,\nabla\ve{u}_{n}\big\|^2 + \frac{\tau\|\nabla\mu\|_{\infty}^2}{\varepsilon\varrho_{\mathrm{min}}\mu_{\mathrm{min}}} \left(\|\hat{\mathbf{U}}_{n}\|^2 + \tau^2\|\chi_{n}\nabla p_{n}\|^2 + \tau\varepsilon\|\sqrt{\mu_n}\,\nabla\ve{u}_{n}\|^2\right)  + \frac{\tau}{\epsilon\mu_{\mathrm{min}}}\|\ve{f}_{n+1}\|_{-1}^2
    \, .
\end{align*}
Adding from $n=0$ to $n=N-1$ and using the Gronwall lemma completes the proof.
}

%\section{Outflow boundary conditions}
%\textcolor{blue}{Comment on (1) natural outflow data and (2) what changes in the Gauge-Uzawa scheme, and why.}

\section{Implementation aspects}\label{sec_impl}
\subsection{Non-Newtonian models}\label{sec_nN}
Variable-density applications often involve generalised Newtonian fluids, where the viscosity may depend on the symmetric velocity gradient---and even on the pressure, as for dense granular flow models \cite{Gesenhues2021}. In that case, using $\nu_{n+1}$ would introduce a non-linearity in the schemes we have presented. This can be easily overcome by using extrapolation, but some care must be taken. Let us take, for example, the BDF2 scheme. An implicit treatment of the viscosity would mean writing
\begin{align*}
    \nu_{n+1} = \nu(\rho_{n+1},|\nabla^{\mathrm{s}}\ve{u}_{n+1}|)\, ,
\end{align*}
where the dependence on $\rho$ determines the phase distribution, and the dependence on $|\nabla^{\mathrm{s}}\ve{u}|$ comes from the rheological model. While it may be tempting to use the interpolation $\nu_{n+1} \approx \nu_{n}^{\star} = 2\nu_n - \nu_{n-1}$, this could yield negative viscosity in some regions. Therefore, we propose instead
\begin{equation*}
    \nu_{n+1} \approx \nu(\rho_{n+1},|\nabla^{\mathrm{s}}\ve{u}_{n}^{\star}|):= \nu^{\star}_{n+1,n}\, ,
\end{equation*}
that is, extrapolating only the velocity field. Of course, we could also extrapolate the density, but there is no need to do that since $\rho_{n+1}$ is already known at the point of tackling the momentum equation. The second part of the viscous term can also be treated in various ways. Two formally second-order extrapolations would be
\begin{align*}
    \nabla^{\top}\ve{u}_{n+1}\nabla\mu_{n+1} \approx 2\nabla^{\top}\ve{u}_{n}\nabla\mu_{n} - \nabla^{\top}\ve{u}_{n-1}\nabla\mu_{n-1} \, ,
\end{align*}
or
\begin{align*}
    \nabla^{\top}\ve{u}_{n+1}\nabla\mu_{n+1} &\approx \nabla^{\top}\ve{u}_{n}^{\star}\nabla(\rho_{n+1}\nu^{\star}_{n+1,n})\, ,
\end{align*}
which is the approach we consider in our simulations. Either way, the stability analysis we derived remains unchanged. On the other hand, there are certain second-order extrapolations, such as $ \nabla^{\top}\ve{u}_{n}^{\star}\nabla\mu^{\star}_{n}$, which may not be as stable because the multiplication of the two extrapolated quantities produces additional, potentially negative terms.

\subsection{Finite element implementation}\label{sec_FEM}
The temporal schemes proposed herein are meant to be general 
and, in principle, usable in combination with various spatial discretisation frameworks. This section discusses some implementation aspects involving finite element discretisations. We will, however, be concise and not dwell on finite element formalism. For the current discussion, it suffices to state that we consider the unknowns $(\ve{u}_{n+1},p_{n+1},\rho_{n+1})$ and their respective test functions to be in suitable, conforming finite element spaces $(X_h,Y_h,Z_h)$, respectively, defined with respect to a triangulation $\mathcal{T}$ of $\Omega$. 

Let us consider the following boundary condition setting for the Navier--Stokes system: 
\begin{align*}
 \ve{u} &= \mathbf{v} && \text{in} \ \ \Gamma_D\times(0,T]\, ,\\ (\mu\nabla\ve{u})\ve{n} - p\ve{n} &= \mathbf{h}  && \text{in} \ \ \Gamma_N\times(0,T]\, , 
\end{align*}
where $(\mathbf{v},\mathbf{h})$ are appropriate data. Taking the first-order coupled scheme as a prototypical setting, the variational formulation for the Navier--Stokes system would be to find $p_{n+1}\in Y_h$ and $\ve{u}_{n+1}\in X_h$, with $\ve{u}_{n+1}|_{\Gamma_D}=\mathbf{v}_{n+1}$, such that
\begin{equation*}
\begin{split}
&\left(\frac{\rho_{n+1}}{\tau}\ve{u}_{n+1} +(\rho_{n+1}\nabla\ve{u}_{n+1})\ve{u}_n + \frac{1}{2}[\nabla\cdot(\rho_{n+1}\ve{u}_n)]\ve{u}_{n+1},\ve{w}\right) + (\mu_{n+1}\nabla\ve{u}_{n+1},\nabla\ve{w}) -(p_{n+1},\nabla\cdot\ve{w})\\
& = \frac{1}{\tau}(\sqrt{\rho_{n}\rho_{n+1}}\,\ve{u}_{n},\ve{w}) 
+\sum_{\Omega_e\in\mathcal{T}}\int_{\Omega_e}(\nabla^{\top}\ve{u}_{n}\nabla\mu_n)\cdot\ve{w}\, \mathrm{d}\Omega + \int_{\Gamma_N}\mathbf{h}_{n+1}\cdot\ve{w}\, \mathrm{d}\Gamma + \langle\ve{f}_{n+1},\ve{w}\rangle \, ,  \\
&(q,\nabla\cdot\ve{u}_{n+1}) = 0\,,
\end{split}
\end{equation*}
for all $q\in Y_h$ and $\ve{w}\in X_h$, with $\ve{w}|_{\Gamma_D}=\ve{0}$. There are two distinct aspects regarding the explicit part of the viscous term. The first is that, differently from the ``Laplacian'' part, the explicit one is \textsl{not} integrated by parts, which is possible because it contains no second-order velocity derivatives. Moreover, we have written the summation over element interiors $\Omega_e$ only since $\mu$ may not have enough regularity for the integral over $\Omega$ to be formally well-defined. That is the case when $\nu$ is a function of $\nabla^{\mathrm{s}}\ve{u}$, which, in the fully discrete case, will be discontinuous (hence $\nabla\mu$ cannot be integrated over internal element boundaries). This formalism, however, has no impact on implementation, since summing over element interiors is exactly how all the other integrals on $\Omega$ are computed. Of course, in the non-Newtonian case, the term containing $\nabla\mu$ is only computable if the velocity is approximated with elements of at least second order. Otherwise, when using linear elements, the viscosity gradient needs to be reconstructed in some other way (e.g., through simple projections or local averaging). 

An attractive feature of our schemes is that, even in the coupled case, the velocity-velocity matrix is block-diagonal. The fully discrete version of the problem above reads 
\begin{align*}
    \begin{bmatrix}
    \mathbf{A}_{n+1} &  \mathbf{B}^{\top} \,    \\
        \mathbf{B}\hphantom{_{n+1}} & \mathbf{0}\hphantom{^{\top}} \,
    \end{bmatrix}
    \begin{bmatrix}
        \ve{U}_{n+1} \\ \ve{P}_{n+1}
    \end{bmatrix} = 
    \begin{bmatrix}
        \ve{F}_{n+1} \\ \ve{0}
    \end{bmatrix},
\end{align*}
where $\mathbf{B}$ is minus the divergence matrix and $\ve{F}_{n+1}$ is the right-hand side vector depending on $(\ve{f}_{n+1}$,$\ve{u}_{n},\mu_{n},\rho_{n+1},\rho_{n},\mathbf{h}_{n+1})$. The velocity-velocity matrix $\mathbf{A}_{n+1}$ consists of $d$ identical blocks:
\begin{align}
\mathbf{A}_{n+1} = 
    \begin{bmatrix}
        \mathbf{K}_{n+1} & \ve{0}\hphantom{_{n+1}}   \\
        \ve{0}\hphantom{_{n+1}} & \mathbf{K}_{n+1}
    \end{bmatrix} 
    \  \text{for} \  d=2\, , \ \, \text{or} \ 
\mathbf{A}_{n+1} = 
    \begin{bmatrix}
        \mathbf{K}_{n+1} & \ve{0}\hphantom{_{n+1}}  & \ve{0}\hphantom{_{n+1}}  \\
         \ve{0}\hphantom{_{n+1}} & \mathbf{K}_{n+1} & \ve{0}\hphantom{_{n+1}} \\
         \ve{0}\hphantom{_{n+1}} & \ve{0}\hphantom{_{n+1}} & \mathbf{K}_{n+1}
    \end{bmatrix} 
    \  \text{for} \  d=3\, ,
    \label{blockDiagonal}
\end{align}
where $\mathbf{K}_{n+1} := \mathbf{M}(\rho_{n+1}) + \mathbf{D}(\mu_{n+1}) + \mathbf{C}(\rho_{n+1},\ve{u}_{n})$, with $\mathbf{M}$, $\mathbf{D}$ and $\mathbf{C}$ denoting, respectively, the mass, diffusion, and convection matrices. The same structure will hold for the BDF2 scheme, with the corresponding extrapolations. The block structure reduces the costs of assembling and solving the linear system \cite{John2016}---especially in the fractional-step case, where it means that the velocity components can be solved separately. 

\section{Numerical examples}\label{sec_examples}
In this section, we assess the accuracy and stability of our methods through various numerical experiments. All the examples use quadrilateral Lagrangian finite elements with first order for pressure and for density, and second order for velocity. In all cases, we write the body force as $\ve{f} = \rho\ve{g}$, where $\ve{g}$ is a given gravity field.

Since the second and third experiments feature very steep density gradients, we use least-squares methods for the density equation, as done by \citet{Pyo2007}. Mind, however, that while least-squares formulations yield unconditional temporal stability when combined with BDF1 or Crank--Nicolson \cite{Pacheco2022CAMWA}, there is currently no theory available on BDF2. 

\subsection{Convergence test}
\begin{figure}[ht!]
\centering
\includegraphics[trim = 0 0 0 0,clip, width = .53\textwidth]{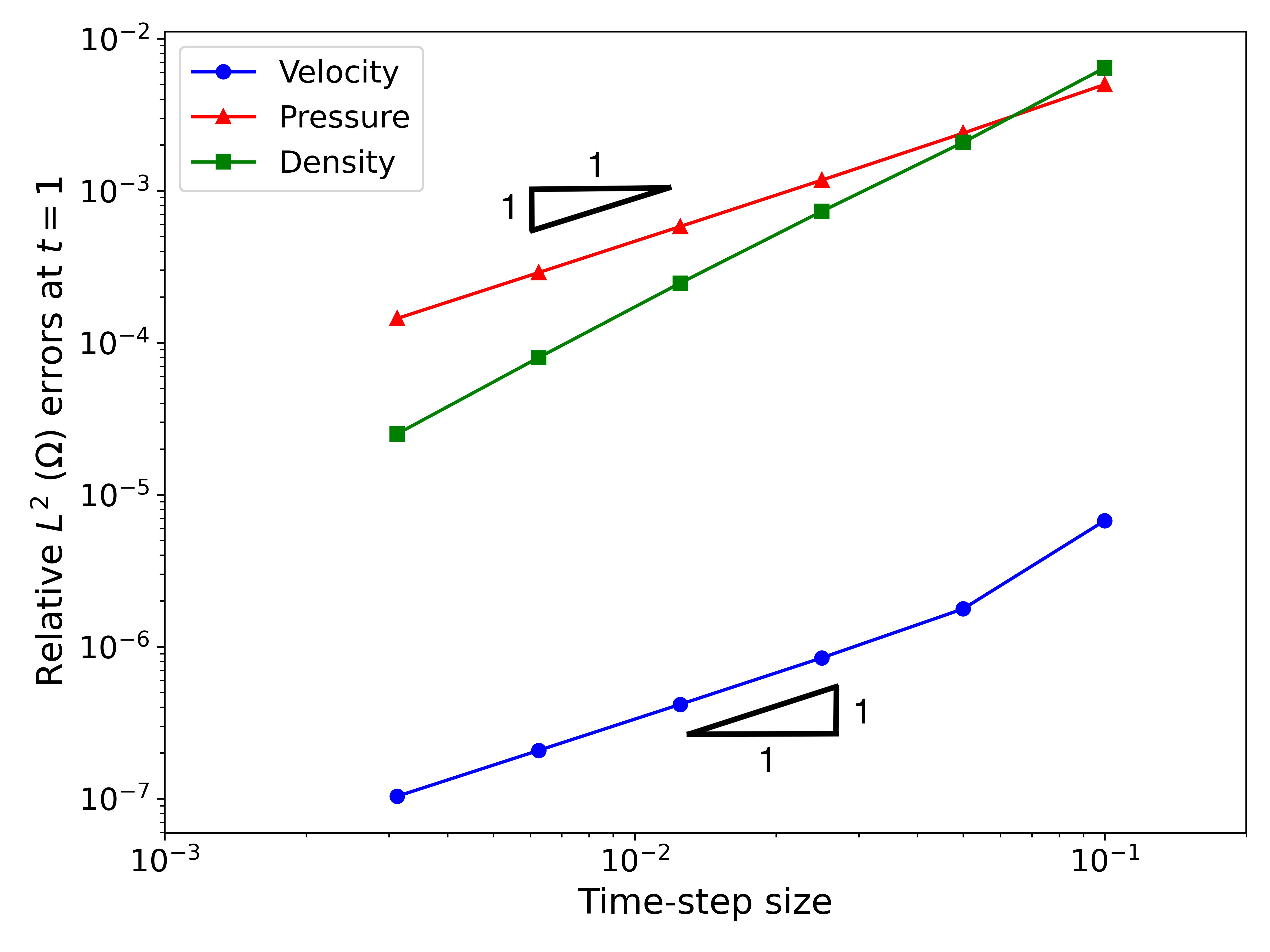}
\caption{Convergence study for the BDF1 scheme, confirming (at least) first-order convergence for all unknowns.}
\label{convergenceBDF1}
\end{figure}

\begin{figure}[ht!]
\centering
\includegraphics[trim = 0 0 0 0,clip, width = .53\textwidth]{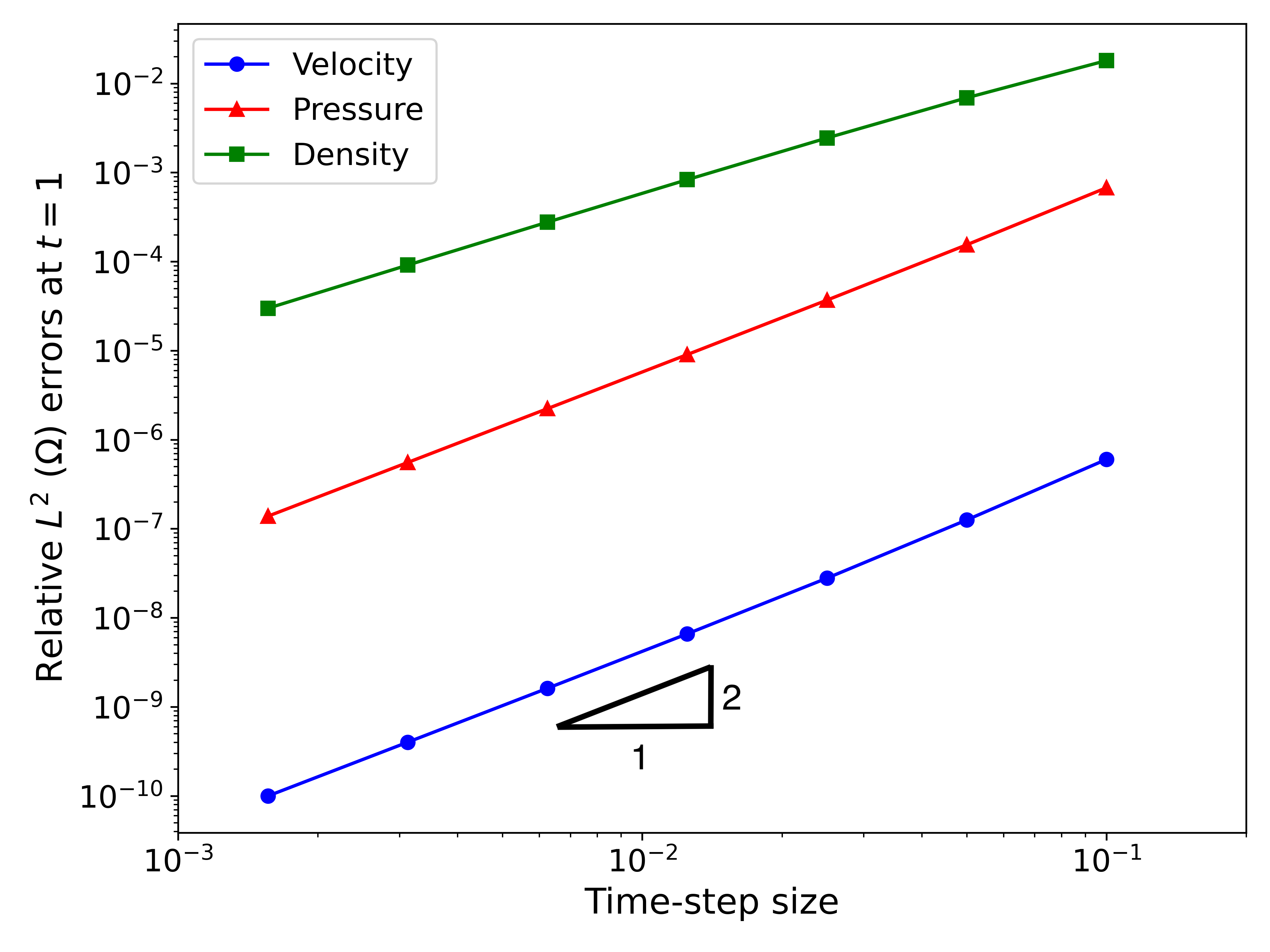}
\caption{Convergence study for the BDF2 scheme, confirming second-order convergence for all unknowns.}
\label{convergenceBDF2}
\end{figure}

\begin{figure}[ht!]
\centering
\includegraphics[trim = 0 0 0 0,clip, width = .53\textwidth]{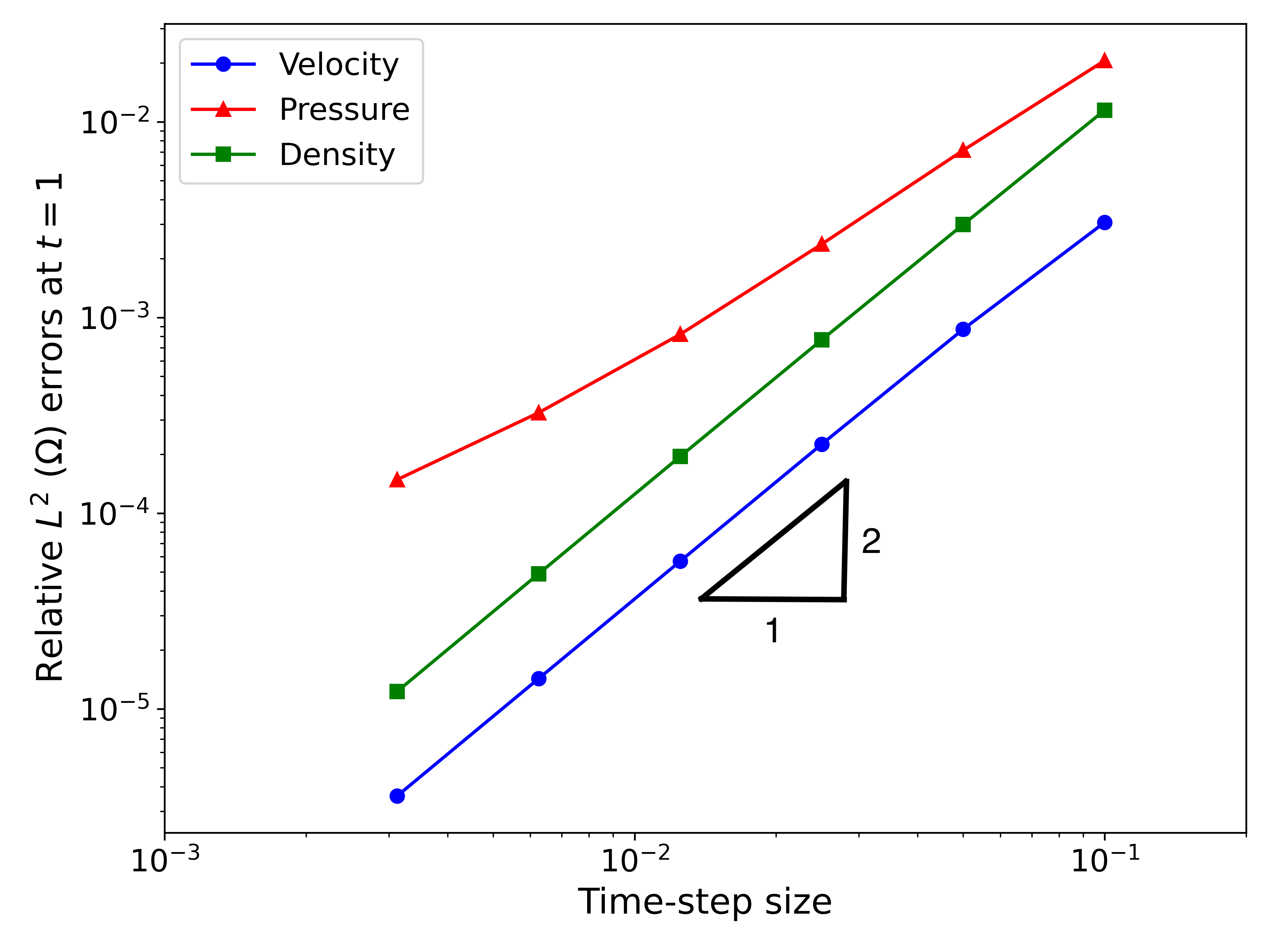}
\caption{Convergence study for the fractional-step scheme, confirming (at least) first-order convergence for all unknowns.}
\label{convergenceBDF1fractional}
\end{figure}

The first problem, solved in $\Omega = (0,1)^2$, has the following analytical solution:
\begin{align*}
\rho &= f(t)\mathrm{e}^{x/u(y)}\,, \\
p &= \sin(2-2x)f(t)\, , \\ 
\ve{u}  &= 
\begin{bmatrix}
u(y)f(t)\\
0
\end{bmatrix}, 
\end{align*}
with the data
\begin{align*}
\ve{g} = 
\begin{bmatrix}
-u(y)f^2(t)\\
0
\end{bmatrix}\, , \ \ \nu = \frac{\mathrm{e}^{-x/u(y)}}{f(t)}\, ,
\end{align*}
where $f(t)=(1+t)^{-1}$ and $u(y) = 1+y-y^2$; the Dirichlet data are computed from the analytical solution. We start with $\tau= 0.1$ and a coarse mesh of $4 \times 4$ square elements and apply five levels of uniform mesh refinement, also halving $\tau$ at each level. As shown in Figures \ref{convergenceBDF1}--\ref{convergenceBDF1fractional}, the convergence rates of all unknowns are quadratic for BDF2 and at least linear for the first-order schemes, as expected. 


\subsection{Rayleigh--Taylor instability}
\begin{figure}[ht!]
\centering
\includegraphics[trim = 0 0 0 0,clip, width = 1\textwidth]{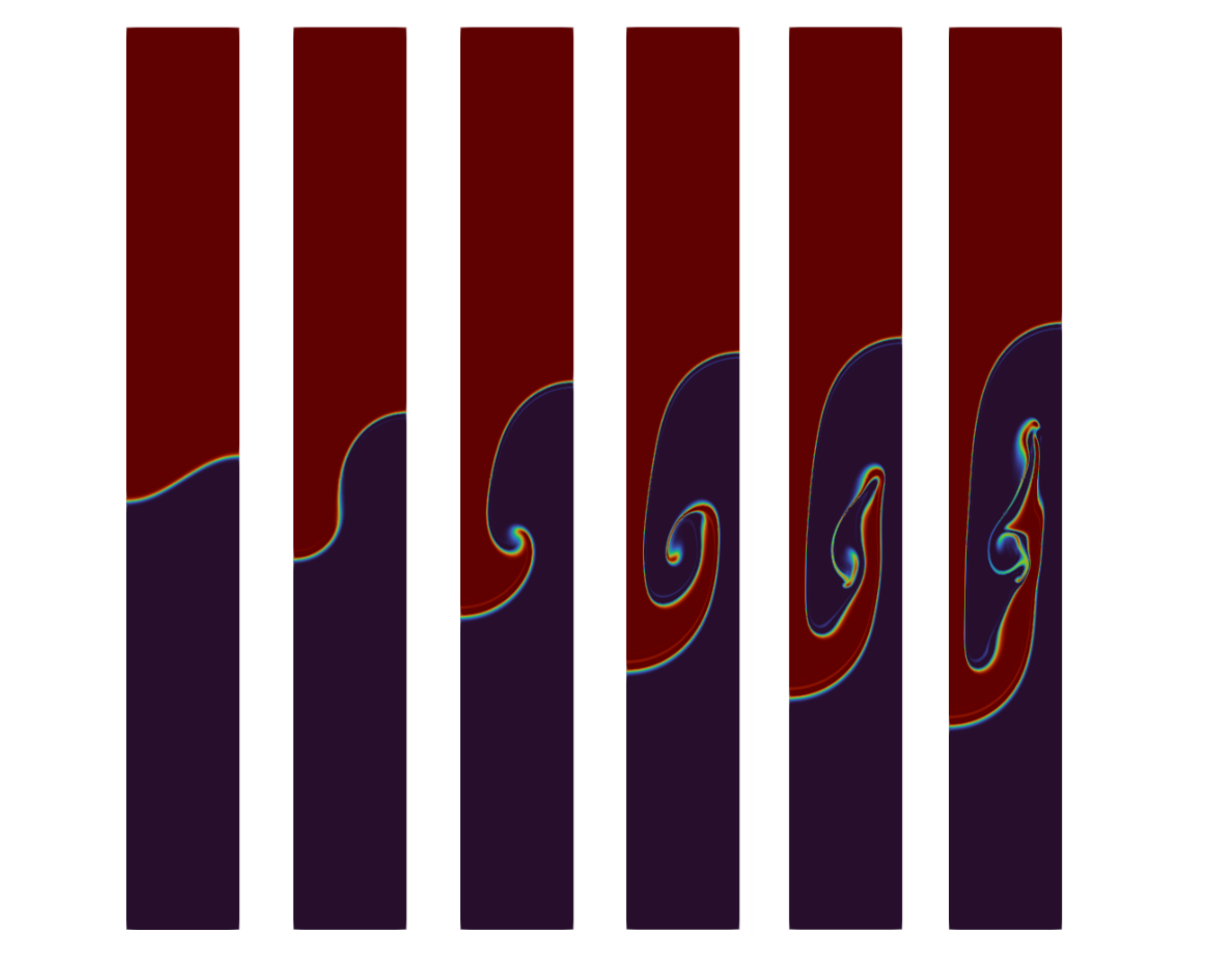}
\caption{Rayleigh--Taylor instability: snapshots of the density field for $t= 0, 1.4, 2.1, 2.45, 2.8$ and $3.15$.}
\label{mushroom}
\end{figure}

The next test is a popular benchmark for variable-density flow solvers, namely the Rayleigh--Taylor instability \cite{Guermond2000}. The setup has two fluids initially at rest in $\Omega = (0,a/2)\times(-2a,2a)$, under gravity $\ve{g} = (0,-g)^{\top}$. Fluid 2 sits initially on top of fluid 1 and is denser ($\rho_2 > \rho_1$). While the standard setup found in the literature uses constant $\mu$, we assume instead constant $\nu$ so that $\mu=\nu\rho$ can vary. The goal is to test the IMEX treatment of variable $\mu$, which is the key novelty of our methods. The initial density field is 
\begin{align*}
\rho|_{t=0} = \frac{\rho_2+\rho_1}{2} + \frac{\rho_2 - \rho_1}{2}\tanh\left(100\frac{y}{a} + 10\cos \frac{2\pi x}{a}\right) .
\end{align*}
 We make the results dimensionless through the following reference quantities: $\nu$ for viscosity, $\rho_1$ for density, $a$ for length, and $\sqrt{ag}$ for velocity, which yields $\sqrt{a/g}$ as temporal scale. The flow regime is fully parametrised by the Reynolds and Atwood numbers:
\begin{align*}
\text{Re} &= \frac{a\sqrt{ag}}{\nu} = 1000\, , \\ \text{At} &= \left|\frac{\rho_2 - \rho_1}{\rho_2 + \rho_1}\right| = \frac{1}{2}\, ,
\end{align*}
which means $\rho_2 = 3\rho_1$. The top and bottom walls are no-slip boundaries, while symmetry (free slip) is enforced on the sides. For the simulation, we use a mesh of $100\times 800$ elements. The time-step size, $\tau=0.001$, is at least an order of magnitude larger than what is often seen in the literature \cite{Guermond2000,Pacheco2022CompMech}, which is possible because our schemes are unconditionally stable. The density field obtained with the second-order scheme is shown in Figure \ref{mushroom} for different times, where we can see some spatial oscillations around the interface; notice, however, that we have not used any type of shock capturing or stabilisation technique.

\begin{figure}[ht!]
\centering
\includegraphics[trim = 0 0 0 0,clip, width = .65\textwidth]{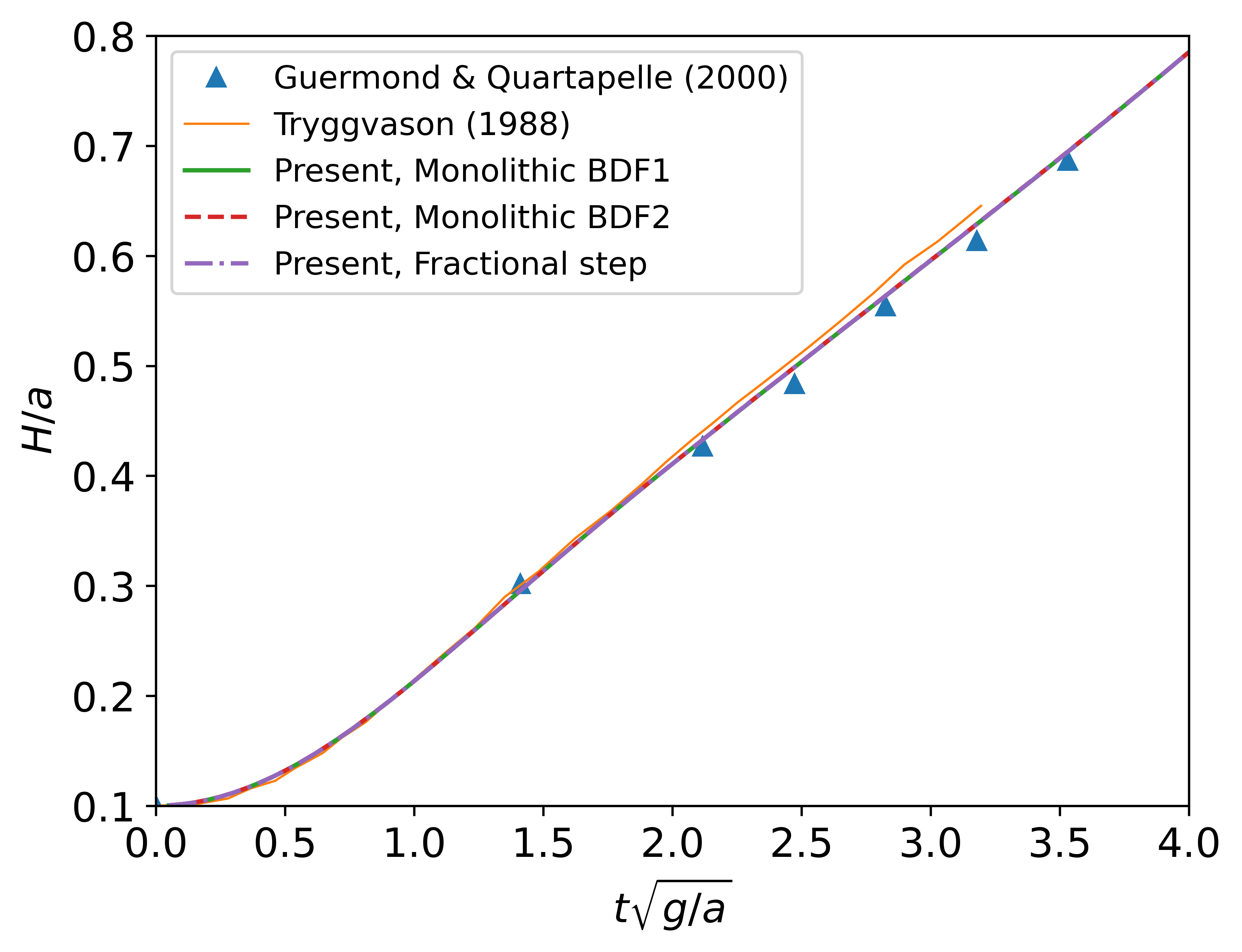}
\caption{Rayleigh--Taylor instability: height of the rising bubble on the right wall.}
\label{mushroomTime}
\end{figure}

Figure \ref{mushroomTime} shows the temporal evolution of the height $H$ of the rising bubble. Our results agree well with the ones reported by \citet{Guermond2000}, even though they have a slightly different setup (constant $\mu$).

\subsection{Viscoplastic falling droplet}\label{sec_droplet}
The last example is inspired by the classical falling droplet benchmark \cite{Calgaro2017}, but we here consider one of the fluids as viscoplastic using the regularised Bingham model \cite{Papanastasiou1987}
\begin{align*}
    \nu(\dot{\gamma}) = \nu_{\infty} + \frac{\sigma_0}{\rho_{\infty}}\frac{1-\mathrm{e}^{-m\dot{\gamma}}}{\dot{\gamma}}\, , \ \ \dot{\gamma} := \sqrt{2\nabla^{\mathrm{s}}\ve{u}:\nabla^{\mathrm{s}}\ve{u}}\, .
\end{align*}
For the non-regularised version ($m\rightarrow\infty$), the material behaves as a Newtonian fluid with viscosity $\nu_{\infty}$ wherever the viscous stress surpasses the yield stress $\sigma_0$, but behaves like a solid ($\nu\rightarrow\infty$) otherwise. Our setup considers a heavy viscoplastic droplet ($\rho_{\infty}=100$) immersed in a lighter fluid ($\rho=1$), which sits on a layer of the viscoplastic fluid (further details on the geometric setup can be found in reference works \cite{Calgaro2017}). The system is initially at rest and then moves under gravitational force $\ve{f} = (0,-\rho)^{\top}$. We set $\nu_{\infty}=10^{-3}$, $\sigma_0 = 1$ and $m=50$. The time-step size is $\tau=0.005$, and the mesh contains $200\times 800$ square elements. 

\begin{figure}[ht!]
\centering
\includegraphics[trim = 0 0 0 0,clip, width = 1\textwidth]{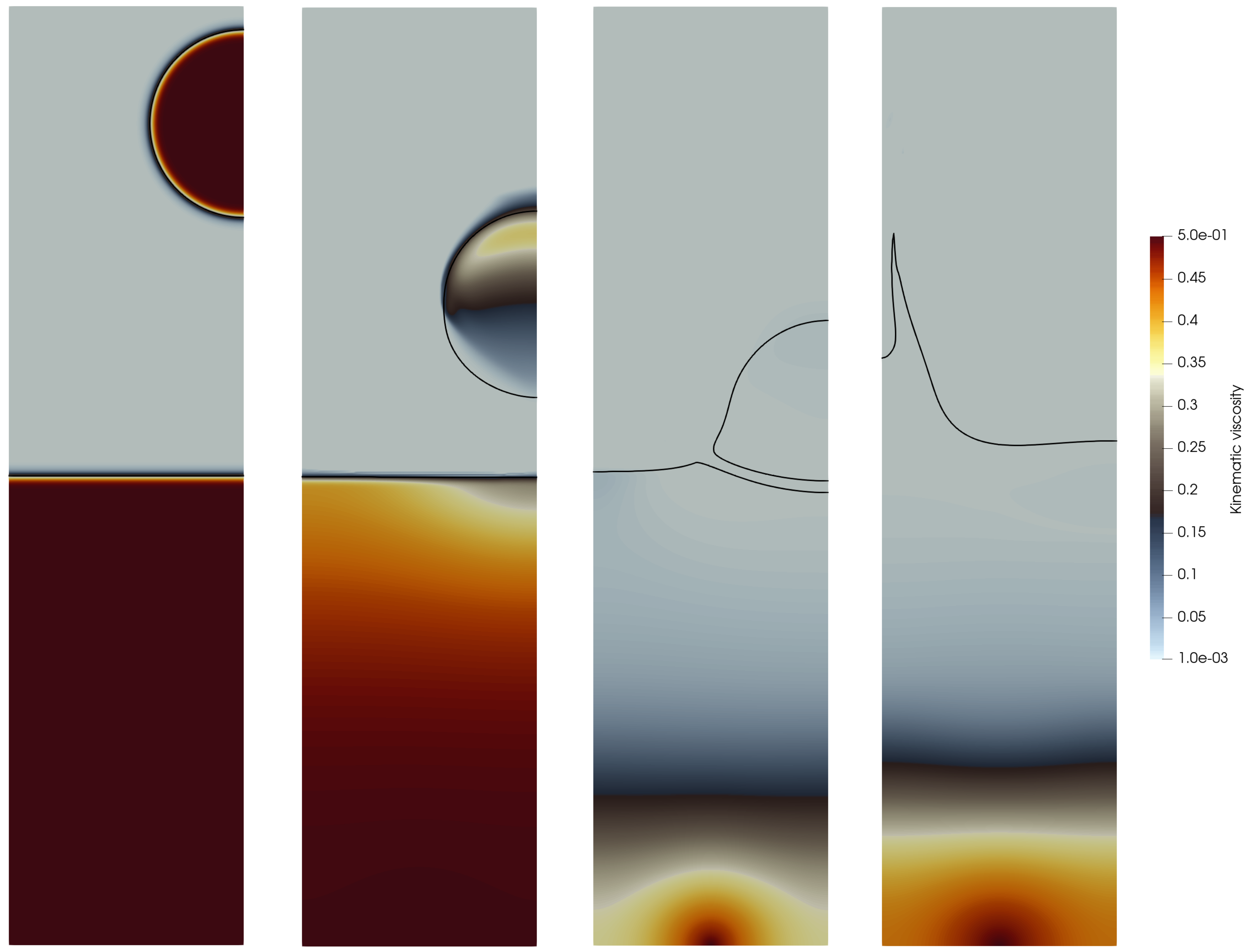}
\caption{Viscoplastic falling droplet: snapshots of the viscosity field $\nu$, for $t= 0, 0.9, 1.15$ and $1.5$. The black contour lines show the interface ($\rho=\frac{100+1}{2}$).}
\label{droplet_snaps}
\end{figure}

This is a very challenging problem to simulate due to the large density ratio between the fluids and the sharp viscosity gradients caused by viscoplasticity. All three IMEX schemes produced similar results, so in Figure \ref{droplet_snaps} we show only the solution of the fractional-step method. The snapshots of the viscosity field reveal that the droplet is essentially solid until shortly before the impact with the lower layer. During the impact, the lower layer fluidises, but then starts hardening again from the bottom. Again, this last example highlights that no numerical instabilities were found and demonstrates that the proposed methodology can cope with complex physical situations.   

\section{Concluding remarks}\label{sec_Conclusion}
In this article, we have devised, analysed and assessed different IMEX schemes for the incompressible variable-density Navier--Stokes system. While most of the literature assumes constant viscosity for the analysis, we have considered the (more realistic) case of variable viscosity, focusing on the IMEX treatment of the diffusive term. Our approach is fully linearised and makes the velocity subsystem block-diagonal, which is especially attractive in the fractional-step case because the resulting system is split, at each time step, into simple scalar subproblems. These features can simplify implementation and increase efficiency. The unconditional temporal stability of our schemes is shown not only through rigorous discrete-in-time analysis but also via challenging numerical examples. Future work will focus on spatial stability, more specifically on positivity-preserving shock-capturing techniques to deal with large density jumps.

\section*{Acknowledgments}
NE and DRQP acknowledge funding by the Federal Ministry of Education and Research (BMBF) and the Ministry of Culture and Science of the German State of North Rhine-Westphalia (MKW) under the Excellence Strategy of the Federal Government and the Länder. The work of GRB has been partially funded by the Leverhulme Trust through the Research Project Grant No. RPG-2021-238.
EC acknowledges the support given by the Agencia Nacional de Investigaci\'on y Desarrollo (ANID) through the project FONDECYT 1210156 and the support given by DICYT-USACH. 

\bibliographystyle{unsrtnat} %order of citation
\bibliography{references}%

\end{document}